\documentstyle[12pt,fleqn,leqno]{article}
\date{April 3,\ 1997}
\title{{\LARGE\bf Sticks and clubs}\\}
\author{Saka\'e Fuchino, Saharon Shelah and Lajos Soukup%
\ifcommented
\bigskip\bigskip\\
\fbox{\parbox{8cm}{\footnotesize\tt\mbox{}\hfill
This is a preliminary version of the paper.\hfill\mbox{}\\
\mbox{}\hfill Any comments are appreciated.\hfill\mbox{}
}}
\fi
}
\newif\iftesting
\newif\ifcommented 
\iftesting
\let\Label\label%
\def\label#1{\mbox{}\marginpar{{\tiny #1}}\Label{#1}}%
\def\xlabel#1{\marginpar{{\tiny #1}}\Label{#1}}%
\fi
\let\xlabel\label
\ifcommented
\fi
{\ifcommented\end{footnotesize}\medskip\\\fi}

\newtheorem{Thm}{{\bf Theorem}}[section]
\newtheorem{Cor}[Thm]{{\bf Corollary}}
\newtheorem{Prop}[Thm]{{\bf Proposition}}
\newtheorem{Lemma}[Thm]{{\bf Lemma}}
\newtheorem{Fact}[Thm]{{\bf Fact}}
\newtheorem{Claim}{{\bf Claim}}[Thm]
\newtheorem{Problem}[Thm]{{\bf Problem}}

\newcommand{\Thmof}[1]{{Theorem \ref{#1}}}
\newcommand{\Corof}[1]{{Corollary \ref{#1}}}
\newcommand{\Propof}[1]{{Proposition \ref{#1}}}
\newcommand{\Lemmaof}[1]{{Lemma \ref{#1}}}
\newcommand{\Factof}[1]{{Fact \ref{#1}}}
\newcommand{\Claimof}[1]{{Claim \ref{#1}}}

\newcommand{\Propabove}{{Proposition \number\theThm}}
\newcommand{\Lemmaabove}{{Lemma \number\theThm}}
\newcommand{\Factabove}{{Fact \number\theThm}}

\newcommand{\prf}{{\bf Proof\ \ }}

\newcommand{\prfofClaim}{\raisebox{-.4ex}{\Large $\vdash$\ \ }}
\newsavebox{\qedbox}\sbox{\qedbox}{
{\unitlength=0.07mm \begin{picture}(40,60)
\put(0,0){\framebox(30,44)[cc]{}}
\put(30,-7){\rule{7\unitlength}{44\unitlength}}
\put(10,-7){\rule{27\unitlength}{7\unitlength}}
\end{picture}}}
\newcommand{\qed}{\mbox{}\hfill\usebox{\qedbox}}
\newcommand{\smallqed}%
{\mbox{}\smallskip\hfill\raisebox{-.4ex}{\Large $\dashv$}\\}
\newcommand{\qedof}[1]%
{\mbox{} \hspace*{\fill}{\usebox{\qedbox}{~(#1)}}%
\mbox{}}
\newcommand{\Qedof}[1]%
{\mbox{} \hspace*{\fill}{\usebox{\qedbox}%
{~(#1~\number\theThm)}}}
\newcommand{\qedofThm}{\Qedof{Theorem}}
\newcommand{\qedofCor}{\Qedof{Corollary}}
\newcommand{\qedofProp}{\Qedof{Proposition}}
\newcommand{\qedofLemma}{\Qedof{Lemma}}
\newcommand{\qedofFact}{\Qedof{Fact}}
\newcommand{\qedskip}{\medskip}
\newcommand{\qedofClaim}%
{\mbox{}\hfill\raisebox{-.4ex}{\Large $\dashv$ }\nolinebreak%
\mbox{~(Claim~\number\theClaim)}}

\def\assert#1{\makebox[5ex][l]{\rm #1)}\ignorespaces}
\def\lassert#1{\llap{\makebox[5ex][l]{\rm #1)}}\ignorespaces}
\def\assertof#1{{\rm #1)}}
\newenvironment{assertion}[1]{\begin{trivlist}
\newbox\assertbox
\dimen255=\textwidth
\setbox\assertbox=\hbox{\hspace*{\parindent}{#1}\hspace{\labelsep}}
\advance\dimen255 by -1\wd\assertbox
\advance\dimen255 by -1ex
\item[]\unhbox\assertbox\hfill
\begin{minipage}[t]{\dimen255}}%
{\end{minipage}\end{trivlist}}
\newcommand{\Implies}{$\,\Rightarrow\,$}
\newcommand{\Equivto}{$\,\Leftrightarrow\,$}
{$\hfill\mbox{}\\[\belowdisplayskip]}
\newcommand{\restr}%
{{\hspace{0.1ex}|\hspace{-0.02ex}{\grave{}}\hspace{0.8ex}}}
\newcommand{\concat}{\mathrel{\mbox{}^\frown}}
\newcommand{\subsetnoneq}%
{\mathrel{\raisebox{-0.8ex}{$\stackrel{\subset}{\scriptstyle\,\not=\,}$}}}
\newcommand{\lessnoneq}%
{\mathrel{\raisebox{-0.8ex}{$\stackrel{<}{\scriptstyle\,\not=\,}$}}}

\newcommand{\powersetof}[1]{{\wp}(#1)}
\newcommand{\cardof}[1]{{\mid{#1}\mid}}

\newcommand{\setof}[2]{\{\,#1\,:\,#2\,\}}
\newcommand{\seqof}[2]{(#1)_{#2}}
\newcommand{\pairof}[1]{(#1)}
\newcommand{\smallsetof}[1]{\{\,#1\,\}}

\newcommand{\continuum}{2^{\aleph_0}}

\newcommand{\forces}[2]{\,\|\hspace{-.35ex}\mbox{\sf--}_{\,#1\,}%
\mbox{\rm``}\,#2\,\mbox{\rm''}}
\newcommand{\notforces}[2]{\rlap{\ /}\|\hspace{-.35ex}\mbox{\sf--}_{\,#1\,}%
\mbox{\rm``}\,#2\,\mbox{\rm''}}
\newcommand{\modelof}[1]{\models\mbox{\rm``}\,#1\,\mbox{\rm''}}
\newcommand{\xmbox}[1]{ ${\rm #1}$ }
\newcommand{\st}{such that}
\newcommand{\wolog}{without loss of generality}
\newcommand{\Wolog}{Without loss of generality}
\newcommand{\wrt}{with respect to}

\newcommand{\tfae}{the following are equivalent}

\newcommand{\cBas}{complete Boolean algebras}

\newcommand{\po}{partial ordering}
\newcommand{\pos}{partial orderings}

\newcommand{\cof}{{\rm cof}}

\newcommand{\Fn}{{\rm Fn}}
\newcommand{\dom}{{\rm dom}}
\newcommand{\otp}{\mathop{\rm otp}}
\newcommand{\supp}{\mathop{\it supp}}
\newcommand{\diff}[1]{\mbox{\it diff\/}(#1)}

\newcommand{\cohenalg}[1]%
{\mbox{$\,\raisebox{0.05ex}{\small$\wr$}\!\!_{_{\!\!}}\mbox{\rm C}_{#1}$}}

\newcommand{\mapping}[3]{#1:#2\rightarrow #3}

\newcommand{\ZFC}{{\rm ZFC}}
\newcommand{\GCH}{{\rm GCH}}
\newcommand{\SCH}{{\rm SCH}}
\newcommand{\CH}{{\rm CH}}
\newcommand{\MA}{{\rm MA}}
\newcommand{\MAC}{{\rm MA(\mbox{\it countable})}}
\newcommand{\MACohen}{{\rm MA(\mbox{\it Cohen}\/)}}

\newcommand{\stickT}{%
\setbox255=\hbox{\raise1ex\hbox{$\hspace{0.2pt}\,\bullet\,$}}
\mathord{\rlap{\hbox to\wd255{\hss\hbox{$|$}\hss}}
\box255}
}
\newcommand{\stickS}{%
\setbox255=\hbox{\raise0.6ex\hbox{$\scriptstyle\bullet$}}
\mathord{\rlap{\hbox to\wd255{\hss\hbox{$\scriptstyle|$}\hss}}
\box255}
}
\newcommand{\stick}{{\mathchoice{\stickT}{\stickT}{\stickS}{\stickS}}}

\newcommand{\calD}{{\cal D}}

\newcommand{\dotb}{{\dot{b}}}

\newcommand{\dotf}{{\dot{f}}}

\newcommand{\dotq}{{\dot{q}}}

\newcommand{\dotu}{{\dot{u}}}
\newcommand{\dotx}{{\dot{x}}}
\newcommand{\doty}{{\dot{y}}}
\newcommand{\dotA}{{\dot{A}}}

\newcommand{\dotG}{{\dot{G}}}
\newcommand{\dotP}{{\dot{P}}}
\newcommand{\dotQ}{{\dot{Q}}}
\newcommand{\dotR}{{\dot{R}}}

\newcommand{\dotU}{{\dot{U}}}

\newcommand{\dotcalP}{\dot{\cal P}}
\begin{document}
\maketitle
\begin{abstract}
We study combinatorial principles known as stick and club. Several  
variants of these principles and cardinal invariants connected to 
them are also considered. 
We introduce a new kind of side-by-side product of \pos\ which we call
pseudo-product. Using such products, we 
give several generic extensions where some of these principles hold 
together with $\neg\CH$ and Martin's Axiom for countable p.o.-sets. 
An iterative version of the pseudo-product is used 
under an inaccessible cardinal to show the consistency of the club 
principle for every stationary subset of limits of $\omega_1$ together 
with $\neg\CH$ and Martin's Axiom for countable p.o.-sets. 
\end{abstract}
{\footnotesize\it
\centerline{Keywords: stick principle, club principle, weak Martin's axiom, 
preservation theorem. }

\centerline{1991 Mathematics Subject classification: $03$ $E$ $35$, 
$03$ $E$ $05$.} }
\section{Beating with sticks and clubs}
In this paper, we study combinatorial principles 
known as `stick' and `club', and their diverse variants 
which are all 
weakenings of $\Diamond$. Hence some of the consequences of 
$\Diamond$ still hold under these principles. 
On the other hand, they are weak enough to be consistent 
with the negation of the continuum hypothesis or even with 
a weak version 
of Martin's axiom in addition.  See e.g.\ 
\cite{broverman-ginsburg-kunen-tall}, \cite{FShS544}, 
\cite{ostaszewski} for applications of these principles. We shall 
begin with introducing the principles and 
some cardinal numbers connected to them. 

$(\stick)$ (read ``stick'') is the following principle introduced 
in S.\ Broverman, J.\ Ginsburg, K.\ Kunen and F.\ 
Tall \cite{broverman-ginsburg-kunen-tall}:
\begin{assertion}{$(\stick)$:}\it
There exists a sequence $(x_\alpha)_{\alpha<\omega_1}$ of 
countable subsets of $\omega_1$ \st\ for any $y\in[\omega_1]^{\aleph_1}$ 
there exists $\alpha<\omega_1$ \st\ $x_\alpha\subseteq y$.
\end{assertion}
Of course the sequence $(x_\alpha)_{\alpha<\omega_1}$ above is a bluff. 
What is essential here is that there exists an 
$X\subseteq[\omega_1]^{\aleph_0}$ of cardinality $\aleph_1$ \st\ for any 
$y\in[\omega_1]^{\aleph_1}$ there is an $x\in X$ with $x\subseteq y$. 
The formulation above is chosen here merely to make the connection to the 
principle $(\clubsuit)$ introduced later, more apparent. 

Note that $(\stick)$ follows from \CH. 
\par
The principle $(\stick)$ suggests the following cardinal number:
\[ \stick= \min
\setof{\cardof{X}}{X\subseteq[\omega_1]^{\aleph_0},\,
\forall y\in[\omega_1]^{\aleph_1}\:\exists x\in X\;x\subseteq y}.
\]\noindent
We have $\aleph_1\leq \stick\leq\continuum$ and $(\stick)$ holds if and 
only if 
$\stick=\aleph_1$. We also consider the following variants of $\stick$:
\[ 
\begin{array}{@{}l@{}l}
\stick'=\min
\setof{\kappa}{&\kappa\geq\aleph_1,\,\mbox{there is an }
X\subseteq[\kappa]^{\aleph_0}\\
&\mbox{\st\ }\cardof{X}=\kappa\mbox{ and }
\forall y\in[\kappa]^{\aleph_1}\:\exists x\in X\;x\subseteq y};
\end{array}
\]\noindent
\[ 
\begin{array}{@{}l@{}l}
\stick''=\min
\setof{\kappa}{&\kappa\geq\aleph_1,\,\mbox{there is an }
X\subseteq[\kappa]^{\aleph_0}\\
&\mbox{\st\ }\cardof{X}=\kappa\mbox{ and }
\forall y\in[\kappa]^{\kappa}\:\exists x\in X\;x\subseteq y};
\end{array}
\]\noindent
\[ 
\begin{array}{@{}l@{}l}
\stick_\lambda=\min
\setof{\cardof{X}}{&X\subseteq[\lambda]^{\aleph_0}\\
&\mbox{\st\ }\forall y\in[\lambda]^{\aleph_1}\:\exists x\in X\;x\subseteq y}.
\end{array}
\]\noindent
We have $\aleph_1\leq \stick''\leq \stick'\leq\continuum$ and 
$\lambda\leq\stick_\lambda\leq\lambda^{\aleph_0}$. 
$(\stick)$ holds if and only if $\stick=\stick'=\stick''=\aleph_1$.
Let us call $X$ as in the definition of $\stick$ 
($\stick'$, $\stick''$ and $\stick_\lambda$ respectively) a 
$\stick$-set ($\stick'$-set, $\stick''$-set and $\stick_\lambda$-set 
respectively). 

\begin{Lemma}\label{st-leq-st'}\mbox{}\\
\noindent
\assert{a} $\stick\leq \stick'$.\\
\assert{b} If $\stick<\aleph_{\omega_1}$ then $\stick=\stick'$.
In particular, we have then $\stick''\leq\stick$. \\
\assert{c} 
If $\lambda\leq\lambda'$ then $\stick_\lambda\leq\stick_{\lambda'}$.\\
\assert{d} $\stick\leq\stick_{\stick}\leq\stick'$.
\end{Lemma}
\prf
\assertof{a}: Let $X\subseteq[\kappa]^{\aleph_0}$ be 
a $\stick'$-set of cardinality $\stick'$. Then 
$X_0=X\cap[\omega_1]^{\aleph_0}$ is a $\stick$-set of cardinality 
$\leq\stick'$. \medskip\\
\assertof{b}: By \assertof{a}, it is enough to show $\stick'\leq\stick$. 
We show inductively that, for every uncountable 
$\kappa\leq \stick$, 
\begin{assertion}{$(*)_\kappa$}\it
there exists an 
$X_\kappa\subseteq[\kappa]^{\aleph_0}$ \st\ 
$\cardof{X_\kappa}\leq\stick$ and \\ 
$\forall y\in[\kappa]^{\aleph_1}\:\exists x\in X_\kappa\;(x\subseteq y)$.
\end{assertion}
For $\kappa=\aleph_1$ this is clear.

Assume that we have shown $(*)_\lambda$ for all $\lambda<\kappa$. If 
$\kappa$ is a successor then by induction
hypothesis, we can 
find $X_\alpha\subseteq[\alpha]^{\aleph_0}$ for all 
$\alpha<\kappa$ \st\ $\cardof{X_\alpha}\leq\stick$ and 
$\forall y\in[\alpha]^{\aleph_1}\:\exists x\in X_\alpha\;x\subseteq y$. 
Let $X_\kappa=\bigcup_{\alpha<\kappa}X_\alpha$. Then $X_\kappa$ has 
the desired property: $\cardof{X_\kappa}\leq\stick$ is clear. If 
$y\in[\kappa]^{\aleph_1}$, there is some 
$\alpha<\kappa$ \st\ $y\in[\alpha]^{\aleph_1}$. Hence 
there is an $x\in X_\alpha\subseteq X$ \st\ $x\subseteq y$.\\
Suppose now that $\kappa$ is a limit. By assumption, we have 
$\cof(\kappa)=\omega$. Let $(\kappa_n)_{n\in\omega}$ be an increasing 
sequence of cardinals below 
$\kappa$ \st\ $\kappa=\bigcup_{n\in\omega}\kappa_n$. For each $n$, let 
$X_{\kappa_n}\subseteq[\kappa_n]^{\aleph_0}$ be as in $(*)_{\kappa_n}$ 
and let 
$X_\kappa=\bigcup_{n\in\omega}X_{\kappa_n}$. Then $X_\kappa$ is as desired: 
clearly $\cardof{X_\kappa}\leq\stick$. If $y\in[\kappa]^{\aleph_1}$ 
there is an $n\in\omega$ \st\ $y\cap\kappa_n$ is uncountable. Hence there 
exists an $x\in X_{\kappa_n}\subseteq X_\kappa$ \st\ 
$x\subseteq y\cap\kappa_n\subseteq y$. 

In particular we have shown that $(*)_\stick$ holds and hence 
$\stick'\leq\stick$. 
\medskip\\
\assertof{c}: Similarly to \assertof{a}.\medskip\\
\assertof{d}: By \assertof{a} and \assertof{c}, we have  
$\stick=\stick_{\aleph_1}\leq\stick_\stick\leq\stick_{\stick'}=\stick'$. 
\qedofLemma\qedskip\\
The question, whether $\stick<\stick'$ is consistent, turned out to be a 
very delicate one: the problem is connected with some natural weakenings 
of \GCH\ whose 
status (i.e.\ whether they are theorems in \ZFC) is still open. One of them  
implies that $\stick=\stick'$ (this is essentially stated in 
\cite[1.2, 1.2A]{Sh:430} in the light of \cite[6.1\,\mbox{[D]}]{Sh:420}; 
for more see \cite{Sh:513}) 
while the negation of the other implies that 
the inequality is consistent. In this paper, we shall treat the latter 
consistency proof 
(\Propof{inequality}). In contrast, the 
consistency of the inequality $\stick''<\stick$ can be shown without any 
such additional set-theoretic assumptions (\Propof{inequality2}). 

The principle $(\clubsuit)$ (`club'), a strengthening of 
$(\stick)$, was 
first formulated in Ostaszewski \cite{ostaszewski}. Let 
$Lim(\omega_1)=\setof{\gamma<\omega_1}{\gamma\mbox{ is a limit}}$. 
For a stationary $E\subseteq Lim(\omega_1)$, 
\begin{assertion}{$\clubsuit(E)$:}\it
There exists a sequence 
$(x_\gamma)_{\gamma\in E}$ of countable 
subsets of $\omega_1$ \st\ for every $\gamma\in E$, $x_\gamma$ is 
a cofinal subset of $\gamma$ with $otp(x_\gamma)=\omega$ and for every 
$y\in[\omega_1]^{\aleph_1}$ there is 
$\gamma\in E$ \st\ $x_\gamma\subseteq y$.
\end{assertion}
Let us call 
$(x_\gamma)_{\gamma\in E}$ 
as above a $\clubsuit(E)$-sequence. 
For $E=Lim(\omega_1)$ we shall simply write $(\clubsuit)$ in place of 
$\clubsuit(Lim(\omega_1))$. 
Clearly $(\stick)$ follows from $(\clubsuit)$. 
Unlike $(\stick)$, $(\clubsuit)$ does not follow from \CH\ since  
$(\clubsuit)$ $+$ \CH\ is known to be equivalent to $\Diamond$ 
(K.\ Devlin, see \cite{ostaszewski}). This equivalence holds also in the 
version argumented with a stationary $E\subseteq Lim(\omega_1)$. 
\begin{Fact}\label{club+ch}
For any stationary $E\subseteq Lim(\omega_1)$, $\clubsuit(E)$ $+$ $\CH$ is 
equivalent to $\Diamond(E)$. 
\end{Fact}
\prf
The proof in \cite{ostaszewski} argumented with $E$ works.
\qedofFact\qedskip

S.\ Shelah \cite{Sh:98} proved the consistency of $\neg\CH$ $+$ 
$(\clubsuit)$ in a model obtained from a model of \GCH\ 
by making the size of $\powersetof{\omega_1}$ to be 
$\aleph_3$ by countable conditions and then collapsing $\aleph_1$ to be 
countable. Soon after that, in an unpublished note, J.\ Baumgartner gave 
a model of $\neg\CH$ $+$ $\clubsuit$ where collapsing of cardinals is 
not involved: his model was obtained from a model of $V=L$ by adding many 
Sacks reals by side by side product. 
I.\ Juh\'asz then proved in an  unpublished note that ``$\neg\CH$ 
$+$ \MAC\ $+$ $(\clubsuit)$'' is consistent. Here \MAC\ 
stands for Martin's 
axiom restricted to countable partial orderings. 
Later P.\ Komj\'ath \cite{komjath} cited a remark by Baumgartner that 
Shelah's model mentioned 
above also satisfies $\neg\CH$ $+$ \MAC\ $+$ $(\clubsuit)$.  
In Section \ref{consis-results}, we shall give yet another model of 
$\neg\CH$ $+$ \MAC\ $+$ 
$(\clubsuit)$ in which collapsing of cardinals is 
not involved (\Thmof{constructible-club}). In section \ref{CSC}, we 
construct a model of $\neg\CH$ $+$ \MAC\ $+$ ``$\clubsuit(E)$ for every 
stationary $E\subseteq Lim(\omega_1)\;$'' starting from a model of \ZFC\ 
with an inaccessible cardinal (\Thmof{clubsuit}). 

These results are rather optimal in the sense that a slight 
strengthening of $\MAC$ implies the negation of 
$(\clubsuit)$. Let $\MACohen$ denote Martin's axiom restricted to the 
partial orderings 
of the form $\Fn(\kappa,2)$ for some $\kappa$ where, as in \cite{kunen}, 
$\Fn(\kappa,2)$ is the p.o.-set for adding $\kappa$ Cohen reals, 
i.e.\ 
the set of functions from some finite subset of $\kappa$ to 2 ordered by 
reverse inclusion. 
\begin{Fact}\label{MACohen} \MA\ for the partial ordering 
$\Fn(\omega_1,2)$ implies $\stick=\stick'=2^{\aleph_0}$. Further, if 
\MACohen\ holds, then we have also $\stick''=2^{\aleph_0}$. 
\end{Fact}
\prf Both equations can be 
proved similarly. For the first equation, it is enough to show 
$\stick=2^{\aleph_0}$ by \Lemmaof{st-leq-st'}. Suppose that 
$X\subseteq[\omega_1]^{\aleph_0}$ is of cardinality less than 
$2^{\aleph_0}$. We 
show that $X$ is not a $\stick$-set. Let $P=\Fn(\omega_1,2)$. Then for 
each $x\in X$ the set 
\[ D_x=\setof{q\in\Fn(\omega_1,2)}{\exists\alpha\in\dom(q)\cap x\ q(\alpha)=0}
\]\noindent
is dense in $P$. For each $\alpha<\omega_1$, 
\[ E_\alpha=
	\setof{q\in\Fn(\omega_1,2)}%
		{\exists\beta>\alpha\,(\beta\in\dom(q)\ \wedge\ q(b)=1)}
\]\noindent
is also a dense subset of $P$. 
Let $\calD=\setof{D_x}{x\in X}\cup\setof{E_\alpha}{\alpha<\omega_1}$ and 
$G$ be a $\calD$-generic filter over $P$. Then the uncountable set 
\[Y=\setof{\alpha<\omega_1}{q(\alpha)=1\mbox{ for some }q\in G}\] 
contains no $x\in X$ as a subset. \qedofFact\qedskip

We shall see in \Propof{inequality2} that \MA\ for the partial ordering 
$\Fn(\omega_1,2)$ is not enough for the last assertion in \Factabove. 

$\clubsuit(E)$ is equivalent to the following seemingly much stronger 
statement. Let $E\subseteq Lim(\omega_1)$ be a stationary set.
\begin{assertion}{$\clubsuit^\dagger(E)$:}\it
 There exists a sequence 
$(x_\gamma)_{\gamma\in E}$ of countable 
subsets of $\omega_1$ \st\ for every $\gamma\in E$, $x_\gamma$ is 
a cofinal subset of $\gamma$ with $otp(x_\gamma)=\omega$ and for every 
$X\in[\omega_1]^{\aleph_1}$,
$\setof{\alpha\in E}{x_\alpha\subseteq X}$ is stationary.
\end{assertion}
\begin{Fact} For any stationary $E\subseteq Lim(\omega_1)$,
$\clubsuit(E)$ and\/ $\clubsuit^\dagger(E)$ are equivalent.
\end{Fact}
\prf
Like \Factof{club+ch}, an easy modification of the corresponding proof in 
\cite{ostaszewski} will work. 
Nevertheless we give here a proof for convenience of the reader. 

Clearly it is enough to show $\clubsuit(E)\Rightarrow\clubsuit^\dagger(E)$. 
Suppose that $(x_\gamma)_{\gamma\in E}$ is a $\clubsuit(E)$-sequence. We 
claim that $(x_\gamma)_{\gamma\in E}$ is then also a 
$\clubsuit^\dagger(E)$-sequence. Otherwise 
there would be a $Y\in[\omega_1]^{\aleph_1}$ and a club 
$C\subseteq Lim(\omega_1)$ 
\st\ $x_\gamma\not\subseteq Y$ for every $\gamma\in C\cap E$. 
By thinning out $C$ if necessary, we may assume that $Y\cap\alpha$ is 
cofinal in $\alpha$ for each 
$\alpha\in C$. For 
$\alpha\in C$, denoting by $\alpha^+$ the next element to $\alpha$ in $C$, 
let 
$y_\alpha\subseteq[\alpha,\alpha^+)\cap Y$ be a 
cofinal subset in $\alpha^+$ with $otp(y_\alpha)=\omega$. Now 
let $Y'=\bigcup_{\alpha\in C}y_\alpha$. Then 
$Y'\in[\omega_1]^{\aleph_1}$ and $Y'\subseteq Y$. We show that 
$\setof{\gamma\in E}{x_\gamma\subseteq Y'}=\emptyset$ which is a 
contradiction: if $\gamma\in E\cap C$ then $x_\gamma\not\subseteq Y'$ follows 
from $Y'\subseteq Y$. If $\gamma\in E\setminus C$ then there is 
$\alpha\in C$ \st\ $\alpha<\gamma<\alpha^+$. By the choice of 
$y_\alpha$, $Y'\cap\gamma$ is not cofinal in $\gamma$. Hence again 
$x_\gamma\not\subseteq Y'$. 
\qedofFact\qedskip\\
\newcommand{\weakclubsuit}{\clubsuit_{\rm w}}
\newcommand{\weakweakclubsuit}{\clubsuit_{\rm w^2}}
\par
Now, let us consider the following variants of the $(\clubsuit)$-principle:
\begin{assertion}{$(\weakclubsuit)$:}\it
 There exists a sequence 
$(x_\gamma)_{\gamma\in Lim(\omega_1)}$ of countable 
subsets of $\omega_1$ \st\ for every $\gamma\in Lim(\omega_1)$, $x_\gamma$ is 
cofinal subset of $\gamma$, $otp(x_\gamma)=\omega$ and for every 
$y\in[\omega_1]^{\aleph_1}$, there is $\gamma<\omega_1$ \st\ 
$x_\gamma\setminus y$ is finite. 
\end{assertion}
\begin{assertion}{$(\weakweakclubsuit)$:}\it
 There exists a sequence 
$(x_\gamma)_{\gamma\in Lim(\omega_1)}$ of countable 
subsets of $\omega_1$ \st\ for every $\gamma\in Lim(\omega_1)$, $x_\gamma$ is 
cofinal subset of $\gamma$, $otp(x_\gamma)=\omega$ and for every 
$y\in[\omega_1]^{\aleph_1}$ 
\[\setof{\alpha<\omega_1}{x_\alpha\cap y\mbox{ is finite}}\cup 
\setof{\alpha<\omega_1}{x_\alpha\setminus y\mbox{ is finite}}\]\noindent  
is stationary in $\omega_1$. 
\end{assertion}
Clearly $(\clubsuit)$ implies $(\weakclubsuit)$. 
Similarly to \Factabove, 
we can prove the equivalence of $(\weakclubsuit)$ with 
$(\weakclubsuit^\dagger)$  
which is 
obtained from $(\weakclubsuit)$ by replacing ``there is an 
$\alpha<\omega_1$ \ldots'' with 
``there are stationary may $\alpha<\omega_1$ \ldots''. Hence 
$(\weakclubsuit)$ implies $(\weakweakclubsuit)$. It is also easy to see 
that $(\weakclubsuit)$ implies $(\stick)$: if 
$(x_\gamma)_{\gamma\in Lim(\omega_1)}$ is a sequence as in the 
definition of $(\weakclubsuit)$, then 
$\setof{x_\gamma\setminus u}{\gamma\in Lim(\omega_1),\,
	u\in[\omega_1]^{<\aleph_0}}$ 
is a $\stick$-set of cardinality $\aleph_1$. 
D\v{z}amonja and Shelah \cite{DjSh} gave a model of $\neg\CH$ $+$ 
$(\weakclubsuit)$ $+$ $\neg(\clubsuit)$. By the remark above this model 
also shows the consistency of non-equivalence of $(\stick)$ and 
$(\clubsuit)$ under $\neg\CH$. 
In this paper we prove that $(\weakweakclubsuit)$ is strictly weaker than 
$(\weakclubsuit)$ by showing the consistency of 
$\neg(\stick)$ $+$ $(\weakweakclubsuit)$ (\Corof{weakweak}). 
The \po\ used in \Corof{weakweak} does not 
force \MAC\ hence the following problem remains open:
\begin{Problem}
Is \MAC\ $+$ $\neg\,(\stick)$ $+$ 
$(\weakweakclubsuit)$ consistent? 
\end{Problem}
\section{Pseudo product of partial orderings}
In this section, we introduce a new kind of side-by-side 
product of p.o.'s which will be used in the next section to prove various 
consistency 
results. Let $X$ be any set and 
$(P_i)_{i\in X}$ be a 
family of \pos. For $p\in\Pi_{i\in X}P_i$ the support of $p$ is defined by 
$\supp(p)=\setof{i\in X}{p(i)\not=1_{P_i}}$. 
For a cardinal $\kappa$, let $\Pi^*_{\kappa,i\in X}P_i$ be 
the set
\[ \setof{p\in\Pi_{i\in X}P_i}{\cardof{\supp(p)}<\kappa}
\]\noindent
with the partial ordering
\[ 
\begin{array}{@{}l@{}l}
p\leq q\quad\Leftrightarrow\quad
	&p(i)\leq q(i)\mbox{ for all }i\in X\mbox{ and}\\
	&\setof{i\in X}{p(i)\lessnoneq q(i)\lessnoneq 1_{P_i}}
		\mbox{ is finite }.
\end{array}
\]\noindent
For $\kappa=\aleph_0$ this is just a finite support product.
We are mainly interested in the case where $\kappa=\aleph_1$.
In this case we shall drop the subscript $\aleph_1$ and write simply 
$\Pi^*_{i\in X}P_i$. 
Further, if $P_i=P$ for some \po\ $P$ for every $x\in X$, we shall write 
$\Pi^*_{\kappa, X}P$ (or even $\Pi^*_{X}P$ when $\kappa=\aleph_1$) to 
denote this \po. 

For $p$, $q\in\Pi^*_{\kappa,i\in X}P_i$ the relation $p\leq q$ can be 
represented as a combination of the two other distinct relations which we 
shall call horizontal and vertical, and denote by $\leq_h$ and $\leq_v$ 
respectively: 
\[ p\leq_h q\quad\Leftrightarrow\quad \supp(p)\supseteq\supp(q)\mbox{ and }
p\restr\supp(q)\subseteq q;
\]\noindent
\[ \begin{array}{@{}l@{}l}
p\leq_v q\quad\Leftrightarrow\quad
	&\supp(p)=\supp(q),\, p(i)\leq q(i)\mbox{ for all }i\in X\mbox{ and}\\
	&\setof{i\in X}{p(i)\lessnoneq q(i)\lessnoneq 1_{P_i}}
		\mbox{ is finite }.
\end{array}
\]\noindent
For $p\in\Pi^*_{\kappa,i\in X}P_i$ and $Y\subseteq X$ let $p\lceil Y$ 
denote the element of $\Pi^*_{\kappa,i\in X}P_i$ defined by 
$p\lceil Y(i)=1_{P_i}$ for every $i\in X\setminus Y$ and 
$p\lceil Y(i)=p(i)$ for $i\in Y$. 

The following is immediate from definition: 
\begin{Lemma}
For $p$, $q\in\Pi^*_{\kappa,i\in X}P_i$, \tfae:\medskip\\
\assert{a} $p\leq q$;\\
\assert{b} There 
is an $r\in\Pi^*_{\kappa,i\in X}P_i$ \st\ $p\leq_h r \leq_v q$;\\
\assert{c}  There 
is an $s\in\Pi^*_{\kappa,i\in X}P_i$ \st\ $p\leq_v s \leq_h q$. 
\end{Lemma}
\prf \assertof{b}\Implies\assertof{a} and 
\assertof{c}\Implies\assertof{a} are clear. For 
\assertof{a}\Implies\assertof{b}, let $r=p\lceil\supp q$; for 
\assertof{a}\Implies\assertof{c}, 
$s=q\restr\supp(q)\cup p\restr(X\setminus\supp(q))$. \qedofLemma
\begin{Lemma}\label{aleph-preserving}\mbox{}\\
\assert{1}
If $P_i$ has the property K for all $i\in X$ then $P=\Pi^*_{i\in X}P_i$ 
preserves $\aleph_1$. \smallskip\\
\assert{2} Suppose that $\lambda\leq\kappa$. 
If $P_i$ has the strong $\lambda$-cc (i.e. 
for every $C\in [P_i]^{\lambda}$ there is pairwise compatible 
$D\in[C]^{\lambda}$), then 
$P=\Pi^*_{\kappa,i\in X}P_i$ preserves $\lambda$.
\end{Lemma}
\prf This proof is a prototype of the arguments we are going to apply
repeatedly.
\assertof{1} and \assertof{2} can be proved similarly. 
For \assertof{1},
assume that there would be $p\in P$ and a $P$-name $\dotf$ \st\medskip\\
\assert{*} 
$p\forces{P}{\mapping{\dotf}{(\omega_1)^V}{\omega}
\mbox{ and }\dotf\mbox{ is 1-1}}$.  \medskip\\
Then, let $(p_\alpha)_{\alpha<\omega_1}$ and 
$(q_\alpha)_{\alpha<\omega_1}$ be 
sequences of elements of $P$ \st\medskip\\
\assert{a} $p_0\leq p$ and $(p_\alpha)_{\alpha<\omega_1}$ is a descending 
sequence \wrt\ $\leq_h$\,;\smallskip\\
\assert{b} $q_\alpha\leq_v p_\alpha$ and 
$q_\alpha$ decides $\dotf(\alpha)$ for all $\alpha<\omega_1$;\smallskip\\
\assert{c} $p_\alpha\restr S_\alpha=q_\alpha\restr S_\alpha$ for every 
$\alpha<\omega_1$ where
\[ S_\alpha=
\supp(q_\alpha)\setminus(\supp(p)\cup\bigcup_{\beta<\alpha}\supp(q_\beta)).
\]\noindent
For $\alpha<\omega_1$ let 
$d_\alpha=\bigcup_{\beta<\alpha}\supp(q_\beta)$. 
Then $(d_\alpha)_{\alpha<\omega_1}$ is a continuously increasing sequence 
in $[X]^{<\omega_1}$. 
Let
$u_\alpha=\setof{\beta\in\supp(q_\alpha)}{q_\alpha(\beta)\not=p_\alpha(\beta)}$ 
for $\alpha<\omega_1$. By \assertof{b}, $u_\alpha$ is finite and by 
\assertof{c} we have  
$u_\alpha\subseteq d_\alpha$. Hence by Fodor's lemma, there exists an 
uncountable (actually even stationary) $Y\subseteq\omega_1$ \st\ 
$u_\alpha=u^*$ for all $\alpha\in Y$, for some fixed 
$u^*\in[X]^{<\aleph_0}$. 
Since $\Pi_{i\in u^*}P_i$ has the property 
K, there exists an uncountable 
$Y'\subseteq Y$ \st\ $\setof{q_\alpha\restr u^*}{\alpha\in Y'}$ 
is pairwise compatible. 
It follows that $q_\alpha$, $\alpha\in Y'$ are pairwise compatible.
For each $\alpha\in Y'$ there exists an 
$n_\alpha\in\omega$ \st\ $q_\alpha\forces{P}{n_\alpha=\dotf(\alpha)}$ 
by \assertof{b}. By \assertof{*}, $n_\alpha$, $\alpha\in Y'$ must be 
pairwise distinct. But this is impossible as $Y'$ is uncountable. 

For \assertof{2}, essentially the same proof works with sequences of 
elements of $P$ of length $\lambda$, using the $\Delta$-system lemma 
argument in place of Fodor's lemma. \\
\qedofLemma
\begin{Lemma}\label{x-cc}
If $\cardof{P_i}\leq 2^{<\kappa}$ for all $i\in X$, then 
$\Pi^*_{\kappa,i\in X}P_i$ has the $(2^{<\kappa})^+$-cc.
\end{Lemma}
\prf By the usual $\Delta$-system lemma argument.\qedofLemma
\begin{Cor}\label{card-preserving}
\assert{a}
Under \CH, if $P_i$ satisfies the property K and 
$\cardof{P_i}\leq\aleph_1$ for every $i\in X$, then $P=\Pi^*_{i\in X}P_i$ 
preserves $\aleph_1$ and has the $\aleph_2$-cc.
In particular $P$ preserves every cardinals.\smallskip\\
\assert{b} Suppose that $2^{<\kappa}=\kappa$. 
If $P_i$ satisfies the strong $\lambda$-cc for every 
$\aleph_1\leq\lambda\leq\kappa$ and $\cardof{P_i}\leq\kappa$ then 
$\Pi^*_{\kappa,i\in X}P_i$ preserves every cardinalities $\leq\kappa$ and has 
the $\kappa^+$-cc. In particular $\Pi^*_{\kappa,i\in X}P_i$ preserves every 
cardinals. 
\end{Cor}
\prf
By Lemmas \ref{aleph-preserving}, \ref{x-cc}. \qedofCor
\begin{Lemma}
For any $Y\subseteq X$ and $x\in X\setminus Y$, we have \label{factorization}
\[{\Pi}^*_{\kappa,i\in X}P_i\cong {\Pi}^*_{\kappa,i\in Y}P_i\times P_x\times
{\Pi}^*_{\kappa,i\in X\setminus(Y\cup\{x\})}P_i.
\] 
\end{Lemma}
\prf
The mapping from ${\Pi}^*_{\kappa,i\in X}P_i$ to 
${\Pi}^*_{\kappa,i\in Y}P_i\times P_x\times
{\Pi}^*_{\kappa,i\in X\setminus(Y\cup\{x\})}P_i$ defined by
\[ p\ \mapsto\ (p\restr Y,p(x),
p\restr(X\setminus(Y\cup\{x\})))
\]\noindent
is an isomorphism. 
\qedofLemma\qedskip\\
In the following we mainly use the \pos\ of the form $\Fn(\lambda,2)$ for 
some 
$\lambda$ as $P_i$ in $\Pi^*_{\kappa,i\in X}P_i$.  
Note that $\Fn(\lambda,2)$ has the property K and strong $\kappa$-cc in the 
sense above for every regular $\kappa$. 

For a pseudo product of the form $\Pi^*_{i\in X}\Fn(\kappa_i,2)$, %
\Lemmaof{aleph-preserving} can be still improved: 
\begin{Thm}{\rm (T.\ Miyamoto)}
For any set $X$, and sequence $(\kappa_i)_{i\in X}$, 
the \po\ $P=\Pi^*_{i\in X}\Fn(\kappa_i,2)$ satisfies the Axiom A.
\end{Thm}
\prf
The sequence of partial orderings $(\leq_n)_{n\in\omega}$ defined by:
$p\leq_0 q\ \Leftrightarrow\ p\leq q$ and $p\leq_n q\ \Leftrightarrow\ %
p\leq_h q\xmbox{\ \  for every }n>0$ witnesses the Axiom A of $P$. 
We omit here the details of the proof since this assertion is never used 
in the following. The idea of the proof needed here is 
to be found in the proof of Lemmas \ref{added-lemma} and \ref{L11}. 
\qedofThm
\begin{Lemma}\label{added-lemma}
Suppose that $\cardof{P_i}\leq\kappa$ for every $i\in X$ and 
$P=\Pi^*_{\kappa^+,i\in X}P_i$. Then\smallskip\\
\assert{1} If $\dotx$ is a $P$-name with $\forces{P}{\dotx\in V}$, then 
for any $p\in P$ there is $q\in P$ \st\ $q\leq_h p$ and
\begin{assertion}{$(\dagger)$} for any $r\leq q$, if $r$ decides 
$\dotx$ then $r\lceil\supp(q)$ already decides $\dotx$.\smallskip\\
\end{assertion}
\assert{2} Let $G$ be $P$-generic. If $u\in V[G]$ is a subset of $V$ of 
cardinality $<\kappa^+$, then there is a ground model set 
$X'\subseteq X$ of cardinality $\leq\kappa$ (in the sense of $V$) \st\ 
$u\in V[G\cap(\Pi^*_{\kappa^+,i\in X'}P_i)]$. 
\end{Lemma}
\prf
\assertof{1}: Let 
$\mapping{\Phi}{\kappa}{\kappa\times\kappa};\,
	\alpha\mapsto (\varphi_1(\alpha),\varphi_2(\alpha))$ 
be a surjection \st\ $\varphi_1(\alpha)\leq\alpha$ for every 
$\alpha<\kappa$. Let $(p_\alpha)_{\alpha<\kappa}$, 
$(p'_\alpha)_{\alpha<\kappa}$ and 
$(r_{\alpha,\beta})_{\alpha<\kappa,\beta<\kappa}$ be sequences of elements 
of $P$ defined inductively by:\medskip\\
\assert{a} $p_0=p$; $(p_\alpha)_{\alpha<\kappa}$ is a descending sequence 
\wrt\ $\leq_h$;\smallskip\\
\assert{b} for a limit $\gamma<\kappa$, $p_\gamma$ is \st\ 
$\supp(p_\gamma)=\bigcup_{\alpha<\gamma}\supp(p_\alpha)$ and, for 
$i\in\supp(p_\gamma)$, $p_\gamma(i)=p_\alpha(i)$ for some 
$\alpha<\gamma$ \st\ $i\in\supp(p_\alpha)$;\smallskip\\
\assert{c} $(r_{\alpha,\beta})_{\beta<\kappa}$ is an enumeration of 
$\setof{r\in P}{r\leq_v p_\alpha}$;\smallskip\\
\assert{d} let $r=r_{\varphi_1(\alpha),\varphi_2(\alpha)}$ and 
\[ p'_\alpha=r\restr\supp(r)\cup p_\alpha\restr(X\setminus\supp(r)).
\]\noindent
If there is $s\leq_h p'_\alpha$ \st\ $s$ decides $\dotx$, then let
\[ p_{\alpha+1}=p_\alpha\restr\supp(p_\alpha)
\cup s\restr(X\setminus\supp(p_\alpha)). 
\]\noindent
Otherwise let $p_{\alpha+1}=p_\alpha$. 

Let $q\in\Pi^*_{\kappa,i\in X}P_i$ be defined by 
$\supp(q)=\bigcup_{\alpha<\kappa}\supp(p_\alpha)$ and, for 
$i\in\supp(q)$, $q(i)=p_\alpha(i)$ for some $\alpha<\kappa$ \st\ 
$i\in\supp(P_\alpha)$. We show that this $q$ is as desired: suppose that 
$r\leq q$ decides $\dotx$. Then there is some $\alpha<\kappa$ \st\ 
\[ r\lceil\supp(q)=p'_\alpha\restr\supp(p'_\alpha)\cup
	q\restr(X\setminus\supp(p'_\alpha)). 
\]\noindent
By \assertof{d}, it follows that 
$r\lceil\supp(q)\leq r\lceil\supp(p_{\alpha+1})$ decides $\dotx$. \medskip\\
\assertof{2}: Let $\dotu$ be a $P$-name for $u$ and let 
$\dotx_\alpha$, $\alpha<\kappa$ be $P$-names \st\ 
$\forces{P}{\dotx_\alpha\in V}$ for every $\alpha<\kappa$ and 
$\forces{P}{u=\setof{\dotx_\alpha}{\alpha<\kappa}}$. By \assertof{1}, for 
each $p\in P$, we can build a sequence $(p_\alpha)_{\alpha<\kappa}$ of 
elements of $P$ decreasing \wrt\ $\leq_h$ \st\ $p_0\leq_h p$ and 
\begin{assertion}{$(\dagger)_{\alpha}$} for any $r\leq p_\alpha$, if 
$r$ decides $\dotx_\alpha$, then $r\lceil\supp(p_\alpha)$ already 
decides $\dotx_\alpha$. 
\end{assertion}
Let $q\in P$ be defined by 
$\supp(q)=\bigcup_{\alpha<\kappa}\supp(p_\alpha)$ and, for 
$i\in\supp(q)$, $q(i)=p_\alpha(i)$ for some $\alpha<\kappa$ \st\ 
$i\in\supp(p_\alpha)$. Then $q$ satisfies:
\begin{assertion}{$(\dagger\dagger)$} for any $r\leq q$, if $r$ decides 
$\dotx_\alpha$ for some $\alpha<\kappa$, then $r\lceil\supp(q)$ already 
decides $\dotx_\alpha$. 
\end{assertion}
The argument above shows that $q$'s with the property 
$(\dagger\dagger)$ are dense in $P$. Hence, by genericity, there is such 
$q\in G$. Clearly, $G\cap\Pi^*_{\kappa,i\in \supp(q)}P_i$ contains every 
information needed to construct $u$. 
\qedofLemma
\section{Consistency results}\xlabel{consis-results}
\begin{Prop}\label{inequality'}
{\rm(\CH)} For any infinite cardinal $\lambda$, let 
$P=\Pi^*_{\lambda}\Fn(\omega_1,2)$. Then 
$\forces{P}{\stick=\lambda}$. 
\end{Prop}
\prf 
\begin{Claim}\label{3.1.1}
$\forces{P}{\stick\geq\lambda}$.\label{stick-geq-lambda}
\end{Claim}
\prfofClaim
If $\lambda=\aleph_1$ this is clear. So assume that $\lambda\geq\aleph_2$. 
For $\xi<\lambda$, let $\dotf_\xi$ be the $P$-name of the generic function 
from $\omega_1$ to $2$ added by the $\xi$-th copy of $\Fn(\omega_1,2)$ in 
$P$. 
Let $G$ be a $P$-generic filter over $V$. In $V[G]$ let 
$X\subseteq[\omega_1]^{\aleph_0}$ be \st\ $\cardof{X}<\lambda$. Then by 
$\aleph_2$-cc of $P$ there exists $\xi<\lambda$ \st\ $X\in V[G']$ for 
$G'=G\cap\Pi^*_{\lambda\setminus\{\xi\}}\Fn(\omega_1,2)$. Since 
$(\dotf_\xi)[G]$ is 
$\Fn(\omega_1,2)$-generic over $V[G']$ by \Lemmaof{factorization}, we have 
$x\not\subseteq((\dotf_\xi)[G])^{-1}\{0\}$ for every $x\in X$.
\qedofClaim
\begin{Claim}\label{3.1.2}
$\forces{P}{\stick\leq\lambda}$.
\end{Claim}
\prfofClaim For $u\in[\lambda]^{<\aleph_0}$, let $\dotcalP_u$ be a 
$P$-name \st\ 
\[ \forces{P}{\dotcalP_u=
([\omega_1]^{\aleph_0})^{V[(\dotf_\xi)_{\xi\in u}]}}
\]\noindent
where $\dotf_\xi$ is as in the proof of the previous claim. Let $\dotcalP$ be 
a $P$-name \st\ 
\[ \forces{P}{\dotcalP=\bigcup\setof{\dotcalP_u}{u\in[\lambda]^{<\aleph_0}}}.
\]\noindent
For each $u\in[\lambda]^{<\aleph_0}$, $(\dotf_\xi[G])_{\xi\in u}$ 
corresponds to a generic filter over 
$\Pi^*_u\Fn(\omega_1,2)\approx \Fn(\omega_1,2)$. 
Hence, by \CH, we 
have $\forces{P}{\cardof{\dotcalP_u}=\aleph_1}$. It follows that 
$\forces{P}{\cardof{\dotcalP}=\lambda}$. 
Thus it is enough to show that 
$\forces{P}{\dotcalP\mbox{ is a }\stick\mbox{-set}}$. 

Let $p\in P$ and $\dotA$ be a 
$P$-name \st\ $p\forces{P}{\dotA\in[\omega_1]^{\aleph_1}}$. 
We show that there is an $r\leq p$ \st\ 
$r\forces{P}{\exists x\in\dotcalP\ x\subseteq\dotA}$. 

Now we proceed as in the proof of \Lemmaof{aleph-preserving}. 
Let $(p_\alpha)_{\alpha<\omega_1}$, $(q_\alpha)_{\alpha<\omega_1}$ be 
sequences of elements of $P$ and  
$(\xi_\alpha)_{\alpha<\omega_1}$ be a strictly increasing sequence of 
ordinals $<\omega_1$ \st\medskip\\
\assert{a} $p_0\leq p$ and $(p_\alpha)_{\alpha<\omega_1}$ is a descending 
sequence \wrt\ $\leq_h$;\smallskip\\
\assert{b} $q_\alpha\leq_v p_\alpha$ and 
$q_\alpha\forces{P}{\xi_\alpha\in\dotA}$ 
for all $\alpha<\omega_1$;\smallskip\\
\assert{c} $p_\alpha\restr S_\alpha=q_\alpha\restr S_\alpha$ for every 
$\alpha<\omega_1$ where
\[ S_\alpha=
\supp(q_\alpha)\setminus(\supp(p)\cup\bigcup_{\beta<\alpha}\supp(q_\beta)).
\]\noindent
For $\alpha<\omega_1$ let 
$u_\alpha=
\setof{\beta\in\supp(q_\alpha)}{q_\alpha(\beta)\not=p_\alpha(\beta)}$. 
As in the proof of \Lemmaof{aleph-preserving}, there 
exists 
$u^*\in[\lambda]^{<\aleph_0}$ \st\ 
$S=\setof{\alpha\in\omega_1}{u_\alpha=u^*}$ is stationary. 
Now $\seqof{q_\alpha\restr u}{\alpha\in S}$ is an infinite sequence of 
elements of 
$P_{u^*}=\Pi_{u^*}\Fn(\omega_1,2)$. Since $P_{u^*}$ satisfies the ccc, there 
exists an $\varepsilon\in S$ and $\zeta<\omega_1$ \st\ 
$q_\varepsilon\restr u^*
\forces{P_{u^*}}{\setof{\xi\in S\cap\zeta}{p_\xi\restr u^*\in \dotG}
\mbox{ is infinite}}$. 
Let 
$r=q_\varepsilon\cup p_\zeta\restr(\supp(p_\zeta)%
\setminus\supp(p_\varepsilon))$. 
Let $\dotb$ be a $P$-name \st\ 
\[ r\forces{P}{\dotb=\setof{\xi\in S\cap\zeta}{q_\xi\restr u^*\in 
\setof{p\restr u^*}{p\in\dotG}}}.
\]\noindent
Let $\dotx$ be 
a 
$P$-name \st\ $r\forces{P}{\dotx=\setof{\xi_\alpha}{\alpha\in\dotb}}$. 
Then $r\forces{P}{\cardof{\dotx}=\aleph_0}$. 
Since $\dotb$ can be computed in $V[(\dotf_\xi[G])_{\xi\in u^*}]$ we have 
$r\forces{P}{\dotx\in\dotcalP_{u^*}}$. It is also clear by definition of 
$\dotx$ that $r\forces{P}{\dotx\subseteq\dotA}$. 
\qedofClaim\\
\qedofProp\qedskip\\
\Propabove\ shows that $\stick$ can be practically every thing. In particular 
we obtain: 
\begin{Cor}
The assertion `$\cof(\stick)=\omega$' is consistent with \ZFC.\qed
\end{Cor}
Actually, $\Fn(\lambda,2)$ forces almost the same situation: 
\begin{Lemma}\label{inequality''}
Suppose that $\lambda$ is a cardinal \st\ $\mu^{\aleph_0}\leq\lambda$ for 
every $\mu<\lambda$. 
Then, for $P=\Fn(\lambda,2)$, we have $\forces{P}{\stick=\lambda}$. 
\end{Lemma}
\prf 
$\forces{P}{\stick\geq\lambda}$ can be proved similarly to 
\Claimof{stick-geq-lambda}. 
For $\forces{P}{\stick\leq\lambda}$, 
let $G$ be a $P$-generic filter and let $G_\alpha=G\cap\Fn(\alpha,2)$ for 
$\alpha<\lambda$. 
In $V[G]$, let
$X=\bigcup\setof{V[G_\alpha]\cap[\omega_1]^{\aleph_0}}{\alpha<\lambda}$. 
Then $\cardof{X}=\lambda$ (here we need \SCH\ in general). We show that 
$X$ is 
a $\stick$-set. For this, it is enough to show the following:
\begin{Claim}
In $V[G]$, if $y\subseteq[\omega_1]^{\aleph_1}$, then there is 
$\alpha^*<\lambda$ and infinite $y'\in V[G_{\alpha^*}]$ \st\ $y'\subseteq y$. 
\end{Claim}
\prfofClaim
In $V$, let $\doty$ be a $P$-name of $y$ which is nice in the sense of 
\cite{kunen}. For $\alpha<\lambda$, let 
$\doty_\alpha=
	\doty\cap\setof{\check{\beta}}{\beta<\omega_1}\times\Fn(\alpha,2)$. 
Then $\forces{P}{\doty=\bigcup_{\alpha<\lambda}\doty_\alpha}$. 
Hence $\forces{P}{\exists\alpha<\lambda\ \doty_\alpha\mbox{ is infinite}}$. 
It follows that there is some $\alpha^*<\lambda$ \st\ 
$y'=\doty_{\alpha^*}[G]$ is 
infinite. 
Since $\doty_{\alpha^*}$ is an $\Fn(\alpha^*,2)$-name, 
$\doty_{\alpha^*}[G]\in V[G_{\alpha^*}]$. Thus these $\alpha^*$ and 
$y'$ are as desired. 
\qedofClaim\\
\qedofLemma
\begin{Prop}\label{inequality}
{\rm (CH)} Suppose that 
\begin{assertion}{$(*)_{\lambda,\mu}$}
There is a sequence 
$(A_i)_{i<\mu}$ of elements of 
$[\lambda]^{\aleph_1}$ \st\ $\cardof{A_i\cap A_j}<\aleph_0$ for every 
$i$, $j<\mu$, $i\not=j$
\end{assertion}
holds for some $\mu>\lambda\geq 2^{\aleph_0}$. 
Then there exists a \po\ $P$ \st\ \medskip\\
\assert{a} $P$ preserves $\aleph_1$ and and has the 
$\aleph_2$-cc; \\
\assert{b} $\forces{P}{\stick=\lambda}$ and \\
\assert{c} $\forces{P}{\stick_\lambda\geq\mu}$.\medskip\\
In particular, 
if $(*)_{\lambda,\mu}$ 
is consistent with \ZFC\ for some $\mu>\lambda\geq 2^{\aleph_0}$, then so 
is $\stick<\stick'$. \end{Prop}
{\bf Remark.\quad}By \cite[\S 6]{Sh:410}, 
$(**)_\mu$ and 
$(*)_{\lambda,\mu}$ for some $\lambda<\mu$ are equivalent, where
\begin{assertion}{$(**)_\mu$} 
there are finite $a_i\subseteq Reg\setminus\aleph_2$ for 
$i<\omega_1$ \st, for any $A\in[\omega_1]^{\aleph_0}$, $\max\, 
pcf(\cup_{i\in A}a_i)\geq\mu$. 
\end{assertion}
For more 
see \cite{Sh:513}. \medskip\\
\prf Let $P$ be as in \Propof{inequality'}. We claim that $P$ is as desired: 
\assertof{a} follows from \Corof{card-preserving} and 
\assertof{b} from \Propof{inequality'}. For \assertof{d}, if 
$X\subseteq[\lambda]^{\aleph_0}$ is a 
$\stick_\lambda$-set then for each $i<\mu$ there is an $x_i\in X$ \st\ 
$x_i\subseteq A_i$. Since $A_i$, $i<\mu$ are almost disjoint $x_i$, 
$i<\mu$ must be pairwise distinct.

The last assertion follows from \Lemmaof{st-leq-st'},\,\assertof{d}. 
\qedofProp \qedskip

Now we show the consistency of the inequality $\stick''<\stick$: 
\begin{Prop}\label{inequality2}
Assume $2^{\aleph_1}=\aleph_2$. Then for any cardinal 
$\lambda\geq\aleph_2$ there exists a \po\ 
$P$ \st\medskip\\
\assert{a} $P$ satisfies the $\aleph_3$-cc;\\
\assert{b} $P$ preserves $\aleph_1$ and $\aleph_2$;\\
\assert{c} if $\lambda^{\aleph_0}=\lambda$ in addition, then 
$\forces{P}{\MA(\Fn(\omega_1,2))}$;\\ 
\assert{d} $\forces{P}{\stick=\lambda}$ and\\
\assert{e} $\forces{P}{\stick''=\aleph_2}$.
\end{Prop}
\prf
\Wolog\ let $\lambda$ be regular and let 
$P=\Pi^*_{\aleph_2,\lambda}\Fn(\omega_1,2)$. Then \assertof{a} 
and \assertof{b} follow from \Corof{card-preserving}. For \assertof{c}, 
note that $\forces{P}{2^{\aleph_0}}=\lambda$ under 
$\lambda^{\aleph_0}=\lambda$. Hence, by \Lemmaof{added-lemma} and 
\Lemmaof{factorization}, we see easily that 
$\forces{P}{\MA(\Fn(\omega_1,2))}$. An argument similar to the 
proof of 
of \Propof{inequality'} shows that 
$\forces{P}{\stick=\lambda}$. For \assertof{e}, we prove first the following: 
\begin{Claim}\label{3.4.1}
Let $X=[\aleph_2]^{\aleph_0}$. Then we have 
$\forces{P}{X\mbox{ is a }\stick''\mbox{-set}}$. In particular 
$\forces{P}{\stick''\leq\aleph_2}$.
\end{Claim}
\prfofClaim
Suppose that, for some $p\in P$ and a 
$P$-name $\doty$ we have $p\forces{P}{\doty\in[\omega_2]^{\aleph_2}}$. 
Let $\dotf$ be a $P$ name \st\ 
$p\forces{P}{\mapping{\dotf}{\omega_2}{\doty}\mbox{ and }\dotf\mbox{ is 1-1}}$.
Let $(p_\alpha)_{\alpha<\omega_2}$ and $(q_\alpha)_{\alpha<\omega_2}$ be 
sequences of elements of $P$ \st\medskip\\
\assert{f} $p_0\leq p$ and $(p_\alpha)_{\alpha<\omega_2}$ is a descending 
sequence \wrt\ $\leq_h$;\smallskip\\
\assert{g} $q_\alpha\leq_v p_\alpha$ and 
$q_\alpha$ decides $\dotf(\alpha)$ for all $\alpha<\omega_2$;\smallskip\\
\assert{h} $p_\alpha\restr S_\alpha=q_\alpha\restr S_\alpha$ for every 
$\alpha<\omega_2$ where
\[ S_\alpha=
\supp(q_\alpha)\setminus(\supp(p)\cup\bigcup_{\beta<\alpha}\supp(q_\beta)).
\]\noindent
For $\alpha<\omega_2$, let $\xi_\alpha\in\omega_2$ be \st\ 
$q_\alpha\forces{P}{\dotf(\alpha)=\xi_\alpha}$. 
Let
$u_\alpha=\setof{\beta\in\supp(q_\alpha)}{q_\alpha(\beta)\not=p_\alpha(\beta)}$ 
for $\alpha<\omega_2$. Just like in the proof of \Lemmaof{card-preserving}, 
we can find 
$u^*\in[\lambda]^{<\aleph_0}$ \st\ $S=\setof{\alpha<\omega_2}{u_\alpha=u^*}$ 
is stationary in $\omega_2$. Since $\cardof{\Fn(\omega_1,2)}=\aleph_1$, 
there exists $T\subseteq S$ of cardinality 
$\aleph_2$ \st\ $q_\alpha\restr u^*$, $\alpha\in T$ are all the same. Let 
$\alpha_n$, $n\in\omega$ be $\omega$ elements of $T$ and let 
$q=\bigcup_{n\in\omega}q_{\alpha_n}$ . Then $q\leq p$ and 
$q\forces{P}{\setof{\xi_{\alpha_n}}{n\in\omega}\subseteq\doty}$. 
\qedofClaim\qedskip\\
Now by \assertof{d}, we have $\forces{P}{\stick''>\aleph_1}$. Hence, by the 
claim above, it follows that $\forces{P}{\stick''=\aleph_2}$.
\qedofProp\qedskip

Modifying the proofs of Propositions \ref{inequality'} and 
\ref{inequality2} slightly, we can also blow up the continuum while setting 
$\stick$ strictly between $\aleph_1$ and $\continuum$. For example:
\begin{Prop}\label{inequality2x}
Assume \CH\ and $2^{\aleph_1}=\aleph_2$. Then for any cardinals 
$\lambda$, $\mu$ \st\ 
$\aleph_2\leq\lambda\leq\mu$ and $\mu^{\aleph_1}=\mu$,  there exists a \po\ 
$P$ \st\medskip\\
\assert{a} $P$ satisfies the $\aleph_3$-cc;\\
\assert{b} $P$ preserves $\aleph_1$ and $\aleph_2$;\\
\assert{c} $\forces{P}{\MAC}$;\\
\assert{d} $\forces{P}{\stick=\lambda}$; \\
\assert{e} $\forces{P}{\stick''=\aleph_2}$ and\\
\assert{f} $\forces{P}{\continuum=\mu}$. 
\end{Prop}
\prf
For $i<\mu$ let 
\[ P_i=\left\{\,
\begin{array}{@{}ll}
\Fn(\omega_1,2), &\mbox{if }i<\lambda,\\
\Fn(\omega,2), &\mbox{otherwise.}
\end{array}
\right.
\]\noindent
Then $P=\Pi^*_{\aleph_2,i<\mu}P_i$ is as desired. \assertof{e} can be proved 
by almost the same proof as that of \Claimof{3.4.1}. \assertof{a}, 
\assertof{b}, \assertof{c} can be shown just as in \Propof{inequality2}. 
Since $P$ adds (at least) $\mu$ many Cohen reals over $V$ and 
$\cardof{P}=\mu$, 
\assertof{f} follows from \assertof{a}. \assertof{d} is proved similarly 
to Claims \ref{3.1.1} and \ref{3.1.2}. For 
$\forces{P}{\stick\leq\lambda}$ we need the following modification of 
\Claimof{3.1.2}: let $\dotcalP$ be defined as in the proof of 
\Claimof{3.1.2}. As there, we can show easily that 
$\forces{P}{\cardof{\dotcalP}=\lambda}$. To show that 
$\forces{P}{\dotcalP\mbox{ is a }\stick\mbox{-set}}$, let $p\in P$ and 
$\dotA$ be a $P$-name \st\ $p\forces{P}{\dotA\in[\omega_1]^{\aleph_1}}$. 
Now let $(p_\alpha)_{\alpha<\omega_1}$, 
$(q_\alpha)_{\alpha<\omega_1}$, $(\xi_\alpha)_{\alpha<\omega_1}$, 
$u^*\in[\mu]^{<\aleph_0}$ and $S$ be just as in the proof of 
\Claimof{3.1.2}. Let $v^*=u^*\setminus\lambda$. Since 
$P_{v^*}=\Pi_{i\in v^*}P_i$ is countable, we may assume \wolog\ that 
$q_\alpha\restr v^*$, $\alpha\in S$ are all the same. Now we can proceed just 
like in the proof of \Claimof{3.1.2} with $u^*$ replaced by  
$u^*\setminus v^*$. \qed\qedskip

The following Lemmas \ref{club-preserving} and \ref{club-forced} show that, 
in spite of typographical similarity, 
$\Pi^*_\lambda\Fn(\omega_1,2)$ and $\Pi^*_\lambda\Fn(\omega,2)$ are 
quite different forcing notions: 
while the first one destroys $(\clubsuit)$ or 
even $(\stick)$ by \Lemmaof{inequality'}, the second one not only preserves 
a $(\clubsuit)$-sequence in 
the ground model but also creates such a sequence generically. 
\begin{Lemma}\label{club-preserving}
Let $S=(x_\gamma)_{\gamma\in E}$ be a 
$\clubsuit(E)$-sequence for a stationary $E\subseteq Lim(\omega_1)$. 
Let $P=\Pi^*_\kappa\Fn(\omega,2)$ for arbitrary $\kappa$. 
Then we have 
$\forces{P}{S\mbox{ is a }\clubsuit(E)\mbox{-sequence}}$.
\end{Lemma}
\prf
Let $p\in P$ and $\dotA$ be a 
$P$-name \st\ $p\forces{P}{\dotA\in[\omega_1]^{\aleph_1}}$. 
We show that there is $q\leq p$ and $\gamma\in E$ 
\st\ $q\forces{P}{x_\gamma\subseteq\dotA}$. 
Let $\dotf$ be a $P$-name \st\ 
$p\forces{P}{\mapping{\dotf}{\omega_1}{\dotA}
\mbox{ and }\dotf\mbox{ is 1-1}}$. Let $(p_\alpha)_{\alpha<\omega_1}$ and 
$(q_\alpha)_{\alpha<\omega_1}$ be sequence of elements of $P$ 
satisfying the conditions \assertof{a} -- \assertof{c} in the proof 
of \Lemmaof{aleph-preserving}. Also, let $u_\alpha$, $\alpha<\omega_1$ be 
as in 
the proof of \Lemmaof{aleph-preserving}. As there, we can find an uncountable 
$Y\subseteq\omega_1$ and $u^*\in[\kappa]^{<\aleph_0}$ \st\ 
$u_\alpha=u^*$ for all $\alpha\in Y$. 
Since $\Pi_{u^*}\Fn(\omega,2)$ is countable we may assume that 
$q_\alpha\restr u^*$ are all the same for $\alpha\in Y$. 
Now for each $\alpha\in Y$ 
let $\beta_\alpha$ be \st\ 
$q_\alpha\forces{P}{\dotf(\alpha)=\beta_\alpha}$ and let 
$Z=\setof{\beta_\alpha}{\alpha\in Y}$. Since $q_\alpha$, $\alpha\in Y$ 
are pairwise compatible, $\beta_\alpha$, $\alpha\in Y$ are pairwise 
distinct and so $Z$ is uncountable. Note that $Z$ is 
a ground model set. Hence there exists $\gamma\in E$ \st\ 
$x_\gamma\subseteq Z$. Let $q=\bigcup_{\alpha\in Y\cap\gamma}q_\alpha$. Then 
$q\leq p$. Since $\sup\setof{\beta_\alpha}{\alpha<\gamma}\geq\gamma$ and 
$\forces{P}{\setof{\beta_\alpha}{\alpha<\gamma}
\xmbox{ is an initial segment of }Z}$, we have 
$q\forces{P}{Z\cap\gamma\subseteq\dotA}$. Hence 
$q\forces{P}{x_\gamma\subseteq\dotA}$. 
\qedofLemma
\begin{Thm}\label{constructible-club}
``$\neg$\CH $+$ \MAC\ $+$ there exists 
a constructible $\clubsuit$\,-sequence'' is consistent.
\end{Thm}
\prf We can obtain a model of the statement by starting from a model of 
$V=L$ and force with $P=\Pi^*_\kappa\Fn(\omega,2)$ for a regular $\kappa$. 
By \Corof{card-preserving}, 
every cardinal of $V$ is preserved in $V[G]$. Since $P$ adds $\kappa$ many 
Cohen reals over $V$ while $\cardof{P}=\kappa$ and $P$ has the 
$\aleph_2$-cc, we have $V[G]\modelof{\continuum=\kappa}$. 
By \Lemmaof{factorization}, $V[G]\modelof{\MAC}$. 
By \Lemmaof{club-preserving}, the $\Diamond$-sequence in $V$ remains a 
$\clubsuit$\,-sequence in $V[G]$.
\qedofThm\qedskip\\

In fact, we do not need a $\clubsuit$-sequence in the ground model to get 
$(\clubsuit)$ in the generic extension by $\Pi^*_{\kappa}\Fn(\omega,2)$\,: 
\begin{Lemma}\label{club-forced}
Let $\kappa$ be uncountable and $P=\Pi^*_\kappa\Fn(\omega,2)$. Then for 
any stationary $E\subseteq Lim(\omega_1)$ we have $\forces{P}{\clubsuit(E)%
\mbox{ holds}}$. 
\end{Lemma}
\prf
For $\gamma\in E$ let 
\[\mapping{f_\gamma}{\,[\gamma,\,\,\gamma+\omega)}%
{\gamma}
\]\noindent 
be a bijection and let 
\[S_\gamma=\setof{x\subseteq\gamma}{x\mbox{ is a cofinal subset of }\gamma,\,
otp(x)=\omega}. 
\]\noindent
For each $x\in S_\gamma$ let $p_x\in P$ be defined by
\[ p_x =
\setof{(\gamma+n,\smallsetof{(0,i)})}%
{n\in\omega,\,i\in 2,\,i=1\,\Leftrightarrow\,
f_\gamma(\gamma+n)\in x}.
\]\noindent
For distinct $x$, $x'\in S_\gamma$, $p_x$ and 
$p_{x'}$ are incompatible. Hence,  
for each $\gamma\in E$, we can find a $ P$-name $\dotx_\gamma$ \st\ 
\[\forces{ P}{\dotx_\gamma\mbox{ is a cofinal subset of }\gamma
\mbox{ and }otp(\dotx_\gamma)=\omega}
\]\noindent 
and 
\[
p_x\forces{ P}{\dotx_\gamma=x}\mbox{\ \ for each }x\in S_\gamma.
\]\noindent 
We show that 
$\forces{P}{(\dotx_\gamma)_{\gamma\in E}
\mbox{ is a }\clubsuit(E)\mbox{-sequence}}$. 
For this, it is enough to show that, for any $p\in P$ and 
a $ P$-name $\dotA$, if 
$p\forces{ P}{\dotA\in [\omega_1]^{\aleph_1}}$, then there is 
$q\leq p$ and $\gamma\in E$ \st\ 
$q\forces{P}{\dotx_\gamma\subseteq\dotA}$. Let $\dotf$
be \st\ 
\[ p\forces{P}{\mapping{\dotf}{\omega_1}{\dotA}\mbox{ and }\dotf
\mbox{ is 1-1}}.
\]\noindent
Now let 
$(p_\alpha)_{\alpha<\omega_1}$, $(q_\alpha)_{\alpha<\omega_1}$, 
$(u_\alpha)_{\alpha\in\omega_1}$, $Y$ and $u^*$ be as in the proof 
of \Lemmaof{aleph-preserving}. For each $\alpha\in Y$ 
let $\beta_\alpha$ be \st\ 
$q_\alpha\forces{ P}{\dotf(\alpha)=\beta_\alpha}$ and let 
$Z=\setof{\beta_\alpha}{\alpha\in Y}$. 
Let 
\[ 
\begin{array}{@{}l@{}l}
C=\setof{\gamma\in Lim(\omega_1)}{\ &%
\bigcup_{\alpha\in Y\cap\gamma}(\supp(q_\alpha)\cap\omega_1)\subseteq\gamma\\
&\mbox{ and }Z\cap\gamma\mbox{ is unbounded in }\gamma}.
\end{array}
\]\noindent
Then $C$ is closed unbounded in $\omega_1$ and hence there exists a 
$\gamma^*\in C\cap E$.
Let $q'=\bigcup_{\alpha\in Y\cap{\gamma^*}}q_\alpha$. Then we have 
$q'\leq q$ and 
$q'\forces{ P}{Z\cap{\gamma^*}\subseteq\dotA}$. Now let 
$x\in S_{\gamma^*}$ be \st\ $x\subseteq Z\cap{\gamma^*}$. Finally let 
$q=q'\cup q_x$. 
Then we have $q\leq p$ and 
$q\forces{ P}{\dotx_\alpha=x\subseteq Z\cap{\gamma^*}\subseteq\dotA}$. 
\qedofLemma
\qedskip\\

Note that $E$'s in Lemmas \ref{club-preserving} and \ref{club-forced} are
ground model sets. To force $\clubsuit(E)$ for every stationary 
$E\subseteq Lim(\omega_1)$ which may be also added generically, we need 
a sort of iteration described in the next section. 

Toward the consistency of $\neg(\weakclubsuit)$ $+$ 
$(\weakweakclubsuit)$, we consider first the following lemma which should 
be a well-known fact. Nevertheless, we include here a proof:
\begin{Lemma}\label{poP}
Assume that there is a sequence 
$(C_\beta)_{\beta<\kappa}$ of elements of 
$[\omega_1]^{\aleph_1}$ \st\ $\cardof{C_\beta\cap C_\gamma}\leq\aleph_0$  
for all $\beta<\gamma<\kappa$. Then there exists a \po\ $P$ with the property 
K \st\ in $V^P$ there is a sequence $(B_\beta)_{\beta<\kappa}$ of elements 
of $[\omega_1]^{\aleph_1}$ \st\ $B_\beta\subseteq C_\beta$ and 
$\cardof{B_\beta\cap B_\gamma}<\aleph_0$ 
for all $\beta<\gamma<\kappa$. 
\end{Lemma}
\prf
Let
\[ \begin{array}{@{}l@{}l}
P = \setof{(D,f)}%
{\ & D\in[\kappa]^{<\aleph_0},\,\mapping{f}{D}{\Fn(\omega_1,2)},\\
 & f(\delta)\in\Fn(C_\delta,2)\mbox{ for all }\delta\in D
}.
   \end{array}
\]\noindent
For $(D,f)$, $(D',f')\in P$, let 
\[ \begin{array}{@{}l@{}l}
(D',f')\leq (D,f)\quad\Leftrightarrow\quad
& D\subseteq D',\,f(\delta)\subseteq f'(\delta)\mbox{ for all }\delta\in D
  \mbox{ and}\\
& (f'(\delta))^{-1}[\{1\}]\setminus (f(\delta))^{-1}[\{1\}],\,\delta\in D
  \mbox{ are pairwise}\\
& \mbox{disjoint}.
   \end{array}
\]\noindent
By the usual $\Delta$-system lemma argument, we can show that $P$ has the 
property K. 
Since $C_\beta$, $\beta<\kappa$ are pairwise disjoint modulo countable, 
the set 
\[ \begin{array}{@{}l@{}l}
\calD_{\beta,\delta}=\setof{(D,f)\in P}%
{
& \beta\in D,\,\delta\in\dom(f(\beta))\mbox{ and}\\
& \exists\eta>\delta\,(\eta\in\dom(f(\beta))\ \wedge\ f(\beta)(\eta)=1)
}
   \end{array}
\]\noindent
is dense in $P$ for every $\beta<\kappa$ and $\delta<\omega_1$. Hence if 
$G$ is a $V$-generic filter over $P$, then
\[ B_\beta=
\setof{\alpha<\omega_1}{f(\beta)(\alpha)=1\mbox{ for some }(D,f)\in G}
\]\noindent
is cofinal in $\omega_1$ and hence uncountable. Also by the definition of 
$\leq$ on $P$, we have 
$\cardof{B_\beta\cap B_\gamma}<\aleph_0$ for every $\beta<\gamma<\kappa$. 
\qedofLemma\qedskip\\
Note that if there is a sequence $(B_\beta)_{\beta<\kappa}$ as 
in \Lemmaabove\ then by the argument in the proof of \Propof{inequality}, 
we have $\stick\geq\kappa$. 
\begin{Lemma}\label{poQ}
There is a \po\ $Q$ with the property K \st\\
$\forces{Q}{(\weakweakclubsuit)}$.
\end{Lemma}
\prf
Let $(Q_\alpha,\dotR_\alpha)_{\alpha\leq\omega_1}$ be the finite support 
iteration of \pos\ with the property K \st\ for each 
$\gamma\in Lim(\omega_1)$, 
there is a 
$Q_\gamma$ name $\dotU_\gamma$ \st\ $Q_\gamma$ forces: 
\begin{assertion}{}
$\dotU_\gamma$ is an ultrafilter over $\gamma$, 
$\gamma\setminus\beta\in\dotU_\gamma$ for all $\beta<\gamma$, 
$\dotR_\alpha$ is a p.o.-set with the property K and there 
is an $\dotR_\gamma$-name $\dotx_\gamma$ \st
\[ \forces{\dotR_\gamma}{\begin{array}[t]{@{}l}
\dotx_\gamma\mbox{ is a cofinal subset of }\gamma
\mbox{ of ordertype }\omega\mbox{ and }\\
\cardof{\dotx_\gamma\setminus a}<\aleph_0\mbox{ for 
all }x\in\dotU_\gamma}.
\end{array}
\]\noindent
\end{assertion}
For example, we can take the Mathias forcing for the 
ultrafilter $\dotU_\gamma$ as  $\dotR_\gamma$. 
For successor $\alpha<\omega_1$ let $\forces{Q_\alpha}{\dotR_\alpha=\{1\}}$. 

Let $Q=Q_{\omega_1}$. 
As $(Q_\alpha,\dotR_\alpha)_{\alpha\leq\omega_1}$ is a finite support 
iteration of property K p.o.s, $Q$ satisfies also the property K (see 
e.g.\ \cite{kunen-tall}). 
Now let $G$ be a $V$-generic filter over $Q$. In 
$V[G]$, if $X\in[\omega_1]^{\aleph_1}$ then the set 
$\setof{\alpha<\omega_1}{X\cap\alpha\in V[G_\alpha]}$ contains a club 
subset $C$ of $Lim(\omega_1)$. Let 
$S_0=\setof{\alpha\in C}{X\cap\alpha\in\dotU_\alpha[G]}$ and 
$S_1=\setof{\alpha\in C}{\alpha\setminus X\in\dotU_\alpha[G]}$. 
Since $\dotU_\alpha[G]$ is an ultrafilter over $\alpha$ in 
$V[G_\alpha]$ for every $\alpha\in C$, we have 
$C=S_0\,\dot\cup\,S_1$. 
We have 
$\cardof{\dotx_\alpha[G]\setminus X}<\aleph_0$ for $\alpha\in S_0$ and 
$\cardof{\dotx_\alpha[G]\cap X}<\aleph_0$ for $\alpha\in S_1$. Thus 
$(\dotx_\alpha[G])_{\alpha\in Lim(\omega_1)}$ is a 
$(\weakweakclubsuit)$-sequence in $V[G]$. 

Actually this proof shows that 
$(\dotx_\alpha[G])_{\alpha\in Lim(\omega_1)}$ is even a 
$(\weakweakclubsuit)$-sequence in the stronger sense that it satisfies 
the assertion of the definition of $(\weakweakclubsuit)$ 
with ``is stationary'' replaced by ``contains a club''. 
\qedofLemma
\begin{Cor}\label{weakweak}
There is a \po\ $R$ with property K \st\ 
$\forces{R}{\stick\geq\aleph_2\mbox{ but }(\weakweakclubsuit)\mbox{ holds}}$. 
In particular $\neg\,(\stick)$ $+$ $(\weakweakclubsuit)$ is consistent 
with \ZFC. 
Further if \CH\ holds then for any cardinal $\kappa$, there exists a 
cardinals preserving proper \po\ 
$R_\kappa$ \st\ 
$\forces{R_\kappa}{\stick\geq\kappa\mbox{ but }
(\weakweakclubsuit)\mbox{ holds}}$. 
\end{Cor}
\prf
Let $R=P_1*\dotP_2$ where $P_1$ is as $P$ in \Lemmaof{poP} for 
$\kappa=\aleph_2$ and $\dotP_2$ as $Q$ in \Lemmaof{poQ} in $V^{P_1}$. 

For the second assertion, we let 
$R_\kappa=\Fn(\kappa,2,\omega_1)*\dotP_1*\dotP_2$. 
Note that under \CH, $\Fn(\kappa,2,\omega_1)$ is cardinals preserving and 
forces that $2^{\aleph_1}=\kappa$. Hence there is a sequence 
$(C_\beta)_{\beta<\kappa}$ as 
in \Lemmaof{poP} in the generic extension. Thus in 
$V^{\Fn(\kappa,2,\omega_1)}$, $\dotP_1$ can be 
taken as in \Lemmaof{poP} for our $\kappa$. Finally, 
in $V^{\Fn(\kappa,2,\omega_1)*\dotP_1}$ 
let $\dotP_2$ be as in \Lemmaof{poQ}. 
\qedofCor
\section{CS$^*$-iteration}
In this section, we introduce an iterative construction of p.o.s which is 
closely related to the pseudo product we introduced in section 2. 
We adopt here 
the conventions of \cite{martin} on forcing. In particular, a p.o.\ 
(or forcing notion) $P$ is a 
{\em pre-ordering} with a greatest element $1_P$. In the following, we 
just try to 
develop a minimal theory needed for \Thmof{clubsuit}. 
More general 
treatment of the iterations like the one described below should be found in 
\cite{properbook}. 

We call a sequence of the form 
$\pairof{P_\alpha,\dotQ_\alpha}_{\alpha\leq\varepsilon}$ a {\it 
CS$^*$-iteration} if the following conditions hold for every 
$\alpha\leq\varepsilon$: 
\begin{assertion}{\phantom{\assertof{1}}}%
\mbox{}\lassert{*0} $P_\alpha$ is a p.o.\ and, if $\alpha<\varepsilon$, 
then $\dotQ_\alpha$ is a $P_\alpha$ name \st\ 
$\forces{P_\alpha}{\dotQ_\alpha
\xmbox{ is a p.o.\ with a greatest element }1_{\dotQ_\alpha}}$. \smallskip\\ 
\mbox{}\lassert{*1} \mbox{}\hfill
$\begin{array}[t]{@{}ll}
        P_\alpha =\setof{p}%
                {&p\mbox{ is a function \st\ }\dom(p)\in[\alpha]^{\leq\aleph_0};\\
                 &p\restr\beta\in P_\beta
                        \mbox{\ \ for any\ \ }\beta<\alpha\mbox{\ \ and,}\\
                 &\mbox{if\ \ }\beta\in\dom(p)\mbox{\ \ then\ \ }
                        p]restr\beta\forces{P_\beta}{p(\beta)\in\dotQ_\beta}}.
   \end{array}
$\hfill\mbox{}\medskip\\
\lassert{*2} For $p$, $q\in P_\alpha$, $p\leq_{P_\alpha} q$ if and only if
\begin{assertion}{\phantom{\assertof{ii}}}%
\mbox{}\lassert{i} for any $\beta<\alpha$, 
$p\restr\beta\forces{P_\beta}{p(\beta)\leq q(\beta)}$;\smallskip\\
\lassert{ii} 
$\diff{p,q}=\setof{\beta\in\dom(p)\cap\dom(q)}%
	{p\restr\beta\notforces{P_\beta}{p(\beta)=q(\beta)}}$ 
is finite.
\end{assertion}
\end{assertion}

We first show that such a sequence 
$\pairof{P_\alpha,\dotQ_\alpha}_{\alpha\leq\varepsilon}$ is really an 
iteration in the usual sense. In the following we assume always that 
$\pairof{P_\alpha,\dotQ_\alpha}_{\alpha\leq\varepsilon}$ is a 
CS$^*$-iteration as defined above. 
\begin{Lemma}\label{L1}
Suppose that $\alpha\leq\beta\leq\varepsilon$. Then\smallskip\\
\assert{0} if $p\in P_\beta$, then $p\restr\alpha\in P_\alpha$;\smallskip\\
\assert{1} $P_\alpha\subseteq P_\beta$;\smallskip\\
\assert{2} for $p$, $q\in P_\alpha$, we have $p\leq_{P_\alpha}q$ \Equivto\ 
$p\leq_{P_\beta}q$;\smallskip\\
\assert{3} for $p$, $q\in P_\beta$, if $p\leq_{P_\beta}q$ then 
$p\restr\alpha\leq_{P_\alpha}q\restr\alpha$. 
\end{Lemma}
\prf \assertof{1} can be proved by induction on $\beta$. Other assertions 
are clear from the definition of CS$^*$-iteration. \qedofLemma
\begin{Lemma}\label{L2}
Suppose that $\alpha\leq\beta\leq\varepsilon$ and $p$, $q\in P_\alpha$. 
Then $p\bot_{P_\alpha}q$ \Equivto\ $p\bot_{P_\beta}q$. 
\end{Lemma}
\prf
Suppose that $p$ and $q$ are compatible in $P_\alpha$, say 
$r\leq_{P_\alpha} p,$ $q$ for some $r\in P_\alpha$. Then $r\in P_\beta$ 
by \Lemmaof{L1},\,\assertof{1} and $r\leq_{P_\beta} p,$ $q$ by 
\Lemmaof{L1},\,\assertof{2}. Hence $p$ and $q$ are compatible in 
$P_\beta$. 

Conversely, suppose that $p$ and $q$ are compatible in 
$P_\beta$, say 
$s\leq_{P_\beta} p,$ $q$ for some $s\in P_\beta$. Then we have 
$s\restr\alpha\in P_\alpha$ by \Lemmaof{L1},\,\assertof{0}, 
$s\restr\alpha\leq_{P_\alpha}p\restr\alpha=p$ and
$s\restr\alpha\leq_{P_\alpha}q\restr\alpha=q$. Hence $p$ and $q$ are 
compatible in $P_\alpha$. 
\qedofLemma\medskip\\
Suppose that $\alpha\leq\beta\leq\varepsilon$, $p\in P_\beta$. 
By \Lemmaof{L1},\,\assertof{0}, we have $p\restr\alpha\in P_\alpha$. For 
$r\leq_{P_\alpha}p\restr\alpha$, let 
\[ p\concat r=p\restr(\dom(p)\setminus\alpha)\,\,\cup\,\, r.
\]\noindent
For $p$, $q\in P_\varepsilon$, 
$p\leq_{P_\varepsilon}^hq$ \Equivto\ $p\leq_{P_\varepsilon}q$ and \label{xxx}
$p\restr\dom(q)=q$; $p\leq_{P_\varepsilon}^vq$ \Equivto\ 
$p\leq_{P_\varepsilon}q$ and $\dom(p)=\dom(q)$ ($h$ and $v$ stand for 
`horizontal' and `vertical' respectively). 
\begin{Lemma}\label{L4} \assert{1} Let $\alpha$, $\beta$, $p$, $r$ be as 
above. Then  $p\concat r\in P_\beta$ and $p\concat r\leq_{P_\beta}r$, 
$p$.\smallskip\\ 
\assert{2} For $p$, $q\in P_\varepsilon$, 
$r=q\restr(\dom(q)\setminus\dom(p))\,\,\cup\,\, p$ is an element of 
$P_\varepsilon$ and $r\leq_{P_\varepsilon}^h p$. \smallskip\\
\assert{3} If $p_n\in P_\varepsilon$ for $n\in\omega$ and 
$p_{n+1}\leq_{P_\varepsilon}^h p_n$ for every $n\in\omega$, then 
$q=\bigcup\setof{p_n}{n\in\omega}$ is an element of $P_\varepsilon$ and 
$q\leq_{P_\varepsilon}^h p_n$ for every $n\in\omega$
\end{Lemma}
\prf
\assertof{1}: By induction on $\beta$. If $\beta=\alpha$ then 
$p\concat r=r\leq p\restr \alpha=p$. Suppose that we have shown the 
inequality for every $\beta'<\beta$. Let $p$ and $r$ be as above. If 
$\beta$ is a limit then we obtain easily $p\concat r\in P_\beta$ and 
$p\concat r\leq_{P_\beta}r$, $p$ by checking \assertof{*1} and 
\assertof{*2} of the definition of CS$^*$-iteration. In particular, 
\assertof{*2},\,\assertof{ii} holds for the inequality 
$p\concat r\leq_{P_\beta}r$, $p$ 
since $\diff{p\concat r,p}=\diff{r,p\restr\alpha}$ and 
$\diff{p\concat r,r}=\emptyset$. 
If $\beta=\gamma+1$ for some $\gamma\geq\alpha$, then 
$p\restr\gamma\concat r\in P_\gamma$, 
$p\restr\gamma\concat r\leq_{P_\gamma}r$, $p\restr\gamma$ by induction 
hypothesis. If $\gamma\not\in\dom(p)$ then it follows 
$p=p\restr\gamma\in P_\beta$ and $p\concat r\leq_{P_\beta}r$, $p$. 
Otherwise 
$(p\concat r)\restr \gamma
        \forces{P_\gamma}{p(\gamma)\leq_{\dotQ_\gamma}p(\gamma)}$. Hence 
again it follows that $p\concat r\in P_\beta$ and 
$p\concat r\leq_{P_\beta}r$, $p$. 

\assertof{2} and \assertof{3} are trivial. 
\qedofLemma
\begin{Lemma}\label{L5} Suppose that $\alpha\leq\beta\leq\varepsilon$, 
$p\in P_\alpha$ and $q\in P_\beta$. If $p$ and $q$ are incompatible in 
$P_\beta$ then $p$ and $q\restr\alpha$ are incompatible in $P_\alpha$. 
\end{Lemma}
\prf Suppose that $p$ and $q\restr\alpha$ are compatible in $P_\alpha$. 
Then there is $r\in P_\alpha$ \st\ $r\leq_{P_\alpha} p$, 
$q\restr\alpha$. Let $s=q\concat r$. By \Lemmaof{L4}, we have 
$s\leq_{P_\beta} q$, $r$. Hence $p$ and $q$ are compatible in $P_\beta$.
\qedofLemma
\begin{Lemma}\label{L6}
Suppose that $\alpha\leq\beta\leq\varepsilon$ and 
that $A$ is a maximal antichain in $P_\alpha$. Then $A$ is also a maximal 
antichain in $P_\beta$. 
\end{Lemma}
\prf
By \Lemmaof{L1},\,\assertof{1}, we have $A\subseteq P_\beta$. By 
\Lemmaof{L2}, $A$ is an antichain in $P_\beta$. Suppose that $A$ were not 
a maximal antichain in $P_\beta$. Then there is some $q\in P_\beta$ \st\ 
$q$ is incompatible with each of $p\in A$. By \Lemmaof{L5}, it follows 
that $q\restr\alpha$ is incompatible with each of $p\restr\alpha=p$, 
$p\in A$. This is a contradiction to the assumption that $A$ is a maximal 
antichain in $P_\alpha$. 
\qedofLemma
\section{CS$^*$-iteration of Cohen reals}\xlabel{CSC}
In the rest, we consider the CS$^*$-iteration 
$\pairof{P_\alpha,\dotQ_\alpha}_{\alpha\leq\kappa}$ for a cardinal 
$\kappa$ \st\ 
\[ \forces{P_\alpha}{\dotQ_\alpha=\Fn(\omega,2)}
\]\noindent
for every $\alpha<\kappa$. 
\begin{Lemma}\label{L10}
Let $p$, $q\in P_\kappa$ be \st\ $p\leq q$. Then there is 
$r\in P_\kappa$ \st\ $r\leq p$ and for any $\alpha\in\diff{r,q}$, there 
is $t\in\Fn(\omega,2)$ \st\ 
$r\restr\alpha\forces{P_\alpha}{r(\alpha)=\check{t}}$. 
\end{Lemma}
\prf We define inductively a decreasing sequence 
$(\alpha_n)_{n<\omega}$ of ordinals and a 
decreasing sequence $(p_n)_{n\in\omega}$ of elements of 
$P_\kappa$ as follows: Let $\alpha_0=\max\diff{p,q}$. Choose 
$p'_0\in P_{\alpha_0}$ so that $p'_0\leq p\restr\alpha_0$ and that 
$p'_0$ decides $p(\alpha_0)$. Let $p_0=p\concat p'_0$. If $\alpha_n$ and 
$p_n$ have been chosen, let $D_n=\diff{p_n,q}\cap\alpha_n$. If 
$D_n=\emptyset$ we are done. Otherwise, 
let $\alpha_{n+1}=\max D_n$. Choose 
$p'_{n+1}\in P_{\alpha_{n+1}}$ \st\ $p'_{n+1}\leq p_n\restr\alpha_{n+1}$ 
and $p'_{n+1}$ decides $p_n(\alpha_{n+1})$. Let 
$p_{n+1}=p_n\concat p'_{n+1}$. This process terminates after $m$ steps 
for some $m\in\omega$, since otherwise we would obtain an infinite 
decreasing sequence of ordinals. Clearly $r=p_m$ is as desired. \qedofLemma

\begin{Lemma}\label{L11} $P_\kappa$ satisfies the axiom A. 
\end{Lemma}
\prf
Let $\leq_n$, $n\in\omega$ be the relations on $P_\kappa$ defined by 
$p\leq_n q$ \Equivto\ $p\leq_{P_\kappa}^h q$ for $p$, 
$q\in P_\kappa$ and {\it every} $n\in\omega$ (in Ishiu 
\cite{ishiu} an axiom A p.o., for which the $\leq_n$'s can be taken to be 
all the same, is called uniformly axiom A).
$\seqof{\leq_n}{n\in\omega}$ has the fusion 
property by \Lemmaof{L4},\,\assertof{3}. Hence it is enough to show the 
following: 
\begin{Claim}
For any $p\in P_\kappa$ and maximal antichain $D\subseteq P_\kappa$, 
there is $q\leq_{P_\kappa}^h p$ \st\ 
$\setof{r\in D}{r\mbox{ is compatible with }q}$ is countable. 
\end{Claim}
\prfofClaim
Let $\mapping{\Phi}{\omega}{\omega\times\omega}$; 
$n\mapsto(\varphi_1(n),\varphi_2(n))$ be a surjection \st\ 
$\varphi_1(n)<n$ for all $n>0$ and, for any $k$, $l\in\omega$, there are 
infinitely many $n\in\omega$ \st\ $\Phi(n)=(k,l)$. We construct 
inductively $p_k$, $t_k$, $u_k\in P_\kappa$ and a sequence 
$\seqof{s_{k,l}}{l\in\omega}$ for 
$k\in\omega$ as follows: let $p_0=p$. If $p_k$ has been chosen then let 
$\seqof{s_{k,l}}{l\in\omega}$ be an enumeration of 
$\Fn(\dom(p_k),\Fn(\omega,2))$. If there are $t\in D$ and 
$u\in P_\kappa$ \st\ $u\leq t$, $p_\kappa$, 
$\diff{u,p_k}=\dom s_{\varphi_1(k),\varphi_2(k)}$ and 
$u\restr\diff{u,p_k}=s_{\varphi_1(k),\varphi_2(k)}$ 
(of course we identify here elements $t$ of $\Fn(\omega,2)$ with 
corresponding $P_\alpha$-name $\check{t}$),
 then let $t_k$ and 
$u_k$ be such $t$ and $u$ and let 
$p_{k+1}=p_k\cup u\restr(\dom(u_k)\setminus\dom(p_k))$. By 
\Lemmaof{L4},\,\assertof{2}, we have $p_{k+1}\in P$. 
Otherwise let 
$t_k=u_k=1_{P_\kappa}$ and $p_{k+1}=k_k$. 

Now, let 
$q=\bigcup_{k\in\omega}p_k$. Then by \Lemmaof{L4},\,\assertof{3}, we have 
$q\in P_\kappa$ and $q\leq_{P_\kappa}p$. We show that this $q$ is as 
desired. 

Suppose that $t\in D$ is compatible with $q$. Then by \Lemmaof{L10}, 
there is $u\subseteq_{P_\kappa}t$, $q$ \st\ $u\restr\diff{q,r}$ has its 
values in $\Fn(\omega, 2)$. Let $n\in\omega$ be \st\ 
$\diff{q,r}\subseteq q_n$ and $k\geq n$ be \st\ 
$s_{\varphi_1(k),\varphi_2(k)}=u\restr\diff{q,r}$. Clearly $t_k\in D$ by 
construction. We claim that $t=t_k$: otherwise $t$ and $t_k$ would be 
incompatible. Hence $u_k$ and $u$ should be incompatible. But this is a 
contradiction. 

It follows that 
\[\setof{r\in D}{r\mbox{ is compatible with }q}\subseteq
        \setof{t_k}{k\in\omega}.\]
\qedofClaim\\
\qedofLemma

In particular, $P_\kappa$ is proper and hence the 
following covering property holds: 
\begin{Cor}\label{CorX}
Suppose that $G$ is a $P_\kappa$-generic filter over $V$. Then for any 
$a\in V[G]$ \st\ 
$V[G]\modelof{a\xmbox{ is a countable set of ordinals}}$, there is a 
$b\in V$ \st\ $a\subseteq b$ and 
$V\modelof{b\xmbox{ is a countable set of ordinals}}$.\qed 
\end{Cor}
\begin{Lemma}\label{L12}
If $\kappa$ is strongly inaccessible, then $P_\kappa$ satisfies the 
$\kappa$-cc. 
\end{Lemma}
\prf Suppose that $p_\beta\in P_\kappa$ for $\beta<\kappa$. We show that 
there are compatible conditions among them. \Wolog\ we may assume that 
$\setof{\dom(p_\beta)}{\beta<\kappa}$ is a $\Delta$-system with the root 
$x\in[\kappa]^{\leq\aleph_0}$ Let 
$\alpha_0=\sup\setof{\gamma+1}{\gamma\in x}$. Then $\alpha_0<\kappa$ and 
$p_\beta\restr x\in P_{\alpha_0}$ for every $\beta<\kappa$. Since 
$\cardof{P_\alpha}<\kappa$ there are $\beta$, $\beta'<\kappa$, 
$\beta\not=\beta'$ \st\ $p_\beta\restr x= p_{\beta'}\restr x$. But then 
$q=p_\beta\cup p_{\beta'}\in P_\kappa$ and $q\leq_{P_\kappa}p_\beta$, 
$p_{\beta'}$. \qedofLemma

\begin{Lemma}\label{L13}
Suppose that $E\subseteq Lim(\omega_1)$ is stationary. Then 
$\forces{P_\kappa}{\clubsuit(E)}$. 
\end{Lemma}
\prf For each $\gamma\in E$ let
$\mapping{f_\gamma}{[\gamma,\gamma+\omega)}{\gamma}$ 
be a bijection and let 
\[ S_\gamma=\setof{x\subseteq\gamma}{x\mbox{ is a cofinal subset of }\gamma,\,
        \otp(x)=\omega}.
\]\noindent
For each $x\in S_\gamma$, let $p_x\in P_\kappa$ be defined by
\[ p_x=\setof{(\gamma+n,\dotq^\gamma_{x,n})}{n\in\omega}
\]\noindent
where $\dotq^\gamma_{x,n}$ is the standard $P_{\gamma+n}$-name for 
$\smallsetof{(0,i)}$ with $i\in 2$ and $i=1$ $\Leftrightarrow$ 
$f_\gamma(\gamma+n)\in x$. For distinct $x$, $x'\in S_{\gamma}$, $p_x$ 
and $p_{x'}$ are incompatible. Hence there is a $P_\kappa$-name 
$\dotx_\gamma$ \st
$\forces{P_\kappa}{\dotx_\gamma\xmbox{ is a cofinal subset of }\gamma
        \xmbox{ with }\otp(\dotx_\gamma)=\omega}$
and 
$p_x\forces{P_\kappa}{\dotx_\gamma=x}$ for every $x\in S_\gamma$. 

We show that 
$\forces{P_\kappa}{\seqof{\dotx_\gamma}{\gamma\in E}\xmbox{ is a }\clubsuit(E)
        \xmbox{-sequence}}$. 
Suppose that $p\in P_\kappa$ and $\dotA$ is a $P_\kappa$-name \st\ 
$p\forces{P_\kappa}{\dotA\in[\omega_1]^{\aleph_1}}$. We have to show that 
there is $q\leq_{P_\kappa}p$ and $\gamma\in E$ \st\ 
$q\forces{P_\kappa}{\dotx_\gamma\subseteq\dotA}$. 

Let $\dotf$ be a $P_\kappa$-name \st\ 
$p\forces{P_\kappa}{\mapping{\dotf}{\omega}{\dotA}\mbox{ is 1-1}}$. 
Choose $p_\alpha$, $q_\alpha$, $u_\alpha$ for $\alpha<\omega_1$ 
inductively \st\ \medskip\\
\assert{a} $p_0\leq_{P_\kappa}p$ and $\seqof{p_\alpha}{\alpha<\omega_1}$ is a 
decreasing sequence \wrt\ $\leq_{P_\kappa}^h$ ;\smallskip\\
\assert{b} $q_\alpha\leq_{P_\alpha}^v p_\alpha$ 
and  $q_\alpha$ decides 
$\dotf(\alpha)$;\smallskip\\ 
\assert{c} 
$u_\alpha=\diff{q_\alpha,p_\alpha}
        \subseteq\dom(p)\cup\bigcup_{\beta<\alpha}\dom(q_\beta)$;\smallskip\\
\assert{d} $q_\alpha\restr{u_\alpha}\in\Fn(\kappa,\Fn(\omega,2))$.\medskip\\
The condition \assertof{d} is possible because of \Lemmaof{L10}. By 
Fodor's lemma, there is $Y\in[\omega_1]^{\aleph_1}$ and 
$r\in\Fn(\kappa,\Fn(\omega,2))$ \st\ $q_\alpha\restr u_\alpha=r$ for 
every $\alpha\in Y$. For each $\alpha\in Y$, there is 
$\beta_\alpha\in\omega_1$ \st\ 
$q_\alpha\forces{P_\kappa}{\dotf(\alpha)=\beta_\alpha}$ by \assertof{b}. 
Let $Z=\setof{\beta_\alpha}{\alpha\in Y}$. 
Let 
\[ 
\begin{array}{@{}l@{}l}
C=\setof{\gamma\in Lim(\omega_1)}{\ &%
\bigcup_{\alpha\in Y\cap\gamma}(\sup(q_\alpha)\cap\omega_1)\subseteq\gamma\\
&\mbox{ and }Z\cap\gamma\mbox{ is unbounded in }\gamma}.
\end{array}
\]\noindent
Then $C$ is closed unbounded in $\omega_1$. Since $E$ was stationary, 
there exists a $\gamma^*\in C\cap E$.
Let $q'=\bigcup_{\alpha\in Y\cap{\gamma^*}}q_\alpha$. Then we have 
$q'\forces{P_\kappa}{Z\cap{\gamma^*}\subseteq\dotA}$. Now let 
$x\in S_{\gamma^*}$ be \st\ $x\subseteq Z\cap{\gamma^*}$. Finally let 
$q=q'\cup p_x$. 
Then we have $q\leq_{P_\kappa}^h p$ and 
$q\forces{P_{\kappa}}{\dotx_\alpha=x\subseteq Z\cap{\gamma^*}\subseteq\dotA}$. 
\qedofLemma
\qedskip

Let 
$\pairof{P_\alpha,\dotQ_\alpha}_{\alpha\leq\kappa}$ be a 
CS$^*$-iteration as above. For $\alpha<\kappa$ let 
$P_\kappa/\dotG_\alpha$ be a $P_\alpha$-name \st\ 
$\forces{P_\alpha}{P_\kappa/\dotG_\alpha=
	\setof{p\in\check{P}_\kappa}{p\restr\alpha\in\dotG_\alpha}
		\xmbox{ with the ordering }p\leq_{\kappa,\alpha}q\ \Leftrightarrow\ 
			p\leq_{P_\alpha}q}$. As in \cite{martin}, we can show that 
$P_\kappa\approx P_\alpha*P_\kappa/\dotG_\alpha$. Also, by \Corof{CorX}, 
practically the same proof as in \cite{martin} shows that 
\[ \forces{P_\alpha}{P_\kappa/\dotG_\alpha\mbox{ is }\approx
	\mbox{ to a CS$^*$-iteration of }\Fn(\omega,2)}.
\]\noindent

Now we are ready to prove the main theorem of this section: 
\begin{Thm}\label{clubsuit}
Suppose that \ZFC\ $+$ ``there exists an inaccessible cardinal'' is 
consistent. Then \ZFC\ $+$ $\neg$\CH\ $+$ \MA$(countable)$ $+$ 
``$\clubsuit(E)$ for every stationary $E\subseteq Lim(\omega_1)$'' is 
consistent as well. 
\end{Thm}
\prf Suppose that $\kappa$ is strongly inaccessible. For $P_\kappa$ as 
above, let $G_\kappa$ be 
a $P_\kappa$-generic filter over $V$. We show that $V[G_\kappa]$ models 
the assertions. 
Let $E\subseteq Lim(\omega_1)$ be a 
stationary set in $V[G_\kappa]$. Since $P_\kappa$ has the $\kappa$-cc by 
\Lemmaof{L12}, there is some $\alpha<\kappa$ \st\ $E\in V[G_\alpha]$ 
where $G_\alpha=G_\kappa\cap P_\alpha$. Hence by the remark before this 
theorem, we may assume \wolog\ that $E\in V$. But then, by \Lemmaof{L13}, 
we have $V[G_\kappa]\modelof{\clubsuit(E)}$. 

Finally, we show that $\MA(countable)$ holds in $V[G_\alpha]$. Let 
$\calD$ be a family of dense subsets of $\Fn(\omega,2)$ in 
$V[G_\kappa]$ of cardinality $<\kappa$. Again by the $\kappa$-cc of 
$P_\kappa$, we can find an $\alpha<\kappa$ \st\ $\calD\in V[G_\alpha]$. 
Since we have
\[ P_\kappa\approx P_\alpha*\dotQ_\alpha*P_\kappa/\dotG_{\alpha+1}\,,
\]\noindent
the generic set over $V[G_\alpha]$ added by 
$\dotQ_\alpha[G_\alpha]=\Fn(\omega,2)$ is $\calD$-generic over 
$\Fn(\omega,2)$ in $V[G_\kappa]$. 
\qedofThm\qedskip

At the moment we --- or more precisely the first and the third author 
--- do not know if an inaccessible cardinal is really necessary in 
\Thmof{clubsuit}. 
As for CS-iteration, $\kappa$ is collapsed to be of cardinality 
$\aleph_2$ in the model 
above, since the continuum of each of the intermediate models is 
collapsed to $\aleph_1$ in the following limit step of cofinality 
$\geq\omega_1$. Thus the following problem seems to be a rather 
hard one: 
\begin{Problem}
Is the combination $\MA(countable)$ $+$ $\clubsuit(E)$ for 
every stationary $E\subseteq Lim(\omega_1)$ consistent with 
$2^{\aleph_0}>\aleph_2$\,? 
\end{Problem}
\mbox{}\bigskip\\
{\large\bf Acknowledgments}

The research of this paper began when the first author (S.F.) was at the 
Hebrew 
University of 
Jerusalem. He would like to thank The Israel Academy of Science and 
Humanities for enabling his stay there. 
He also would like to thank T.\ Miyamoto for some quite helpful 
remarks. 

The second author (S.S.) was partially supported by the Deutsche 
Forschungs\-gemein\-schaft(DFG) grant Ko 490/7--1. 
He also gratefully acknowledges partial support by the Edmund Landau Center 
for research in Mathematical Analysis, supported by the Minerva 
Foundation (Germany). 
The present paper is the second author's 
Publication No.\ 544. 

The third author (L.S.) is partially supported by the Hungarian National 
Foundation for Scientific Research grant No.\ 16391 and the Deutsche 
Forschungs\-gemein\-schaft (DFG) grant Ko 490/7--1.


\newpage
\mbox{}\vfill\mbox{}\hfill
\parbox[t]{7cm}{
\noindent
{\bf Authors' addresses}\bigskip\bigskip\\
\it
{\rm Saka\'e Fuchino}\smallskip\\
Institut f\"ur Mathematik II,\\ Freie Universit\"at Berlin\\
14195 Berlin, Germany\medskip\\
{\tt fuchino@math.fu-berlin.de}
\bigskip\bigskip\\
{\rm Saharon Shelah}\smallskip\\
Institute of Mathematics,\\ The Hebrew University of Jerusalem\\
91904 Jerusalem, Israel\smallskip\\
and\smallskip\\
Department of Mathematics,\\ Rutgers University\\
New Brunswick, NJ 08854, USA\medskip\\
{\tt shelah@math.huji.ac.il}
\bigskip\bigskip\\
{\rm Lajos Soukup}\smallskip\\
Mathematical Institute\\ of the Hungarian Academy of Sciences\medskip\\
{\tt soukup@math-inst.hu}}

\end{document}